\documentclass[leqno,12pt]{amsart}
\usepackage{amssymb} 
\theoremstyle{plain} 
\newtheorem{theorem}{\indent\sc Theorem}[section] 

\newtheorem{corollary}[theorem]{\indent\sc Corollary}
\newtheorem{proposition}[theorem]{\indent\sc Proposition}

\newtheorem{question}[theorem]{\indent\sc Question}
\theoremstyle{definition} 

\newtheorem{remark}[theorem]{\indent\sc Remark}

\begin{document}

\title[Symmetric Hamiltonian of the Garnier system]{Symmetric Hamiltonian of the Garnier system and its degenerate systems in two variables \\}

\author{By\\
Yusuke Sasano}

\renewcommand{\thefootnote}{\fnsymbol{footnote}}
\footnote[0]{2000\textit{ Mathematics Subjet Classification}.
34M55; 34M45; 58F05; 32S65.}

\maketitle

\begin{abstract} We present {\it symmetric Hamiltonians} for the degenerate Garnier systems in two variables. For these symmetric Hamiltonians, we make the symmetry and holomorphy conditions, and we also make a generalization of these systems involving symmetry and holomorphy conditions inductively. We also show the confluence process among each system by taking the coupling confluence process of the Painlev\'e systems.
\end{abstract}

\section{Introduction}

In this paper, we consider the Garnier system G(1,1,1,1,1) and its degenerate systems in two variables, which are completely integrable Hamiltonian systems of the form
\begin{align}\label{into1}
\begin{split}
&dx=\frac{\partial H_1}{\partial y}dt+\frac{\partial H_2}{\partial y}ds, \quad dy=-\frac{\partial H_1}{\partial x}dt-\frac{\partial H_2}{\partial x}ds,\\
&dz=\frac{\partial H_1}{\partial w}dt+\frac{\partial H_2}{\partial w}ds, \quad dw=-\frac{\partial H_1}{\partial z}dt-\frac{\partial H_2}{\partial z}ds.
\end{split}
\end{align}
The Hamiltonians $H_1,H_2$ are polynomial with respect to $x,y,z,w$ whose coefficients are rational functions of $t,s$.

As is explained in \cite{K2}, these systems are obtained as monodromy preserving deformations equations of second-order linear ordinary differential equations with regular or irregular singular points and apparent singular points. Let us assign 1 to a regular singular point and r+1 to an irregular singular point of Poincar\'e rank r.

In this paper, we present the polynomial Hamiltonian system \eqref{into1} with {\it symmetric Hamiltonians}
\begin{equation}
H_1:=H_{*}(x,y,t;\alpha_0,\ldots)+R(x,y,z,w,t,s;\alpha_0,\ldots) \in {\Bbb C}(t,s)[x,y,z,w], \quad H_2=\pi(H_1),
\end{equation}
where the transformation $\pi$ is given by
\begin{equation}
\pi:(x,y,z,w,t,s) \rightarrow (z,w,x,y,s,t), \quad (\pi)^2=1
\end{equation}
with some parameter's change, and the symbol $H_{*}(x,y,t;\alpha_0,\ldots)$ denotes one of the Painlev\'e Hamiltonians. The birational symplectic transformations take known Hamiltonian systems (see \cite{K2}) to the symmetric Hamiltonian systems. These symmetric Hamiltonians are new.

For the symmetric Hamiltonian systems, we make the symmetry and holomorphy conditions.

As is well-known, the degeneration from $P_{VI}$ to $P_{V}$ (see \cite{T2}) is given by
\begin{gather*}
\begin{gathered}
\alpha_0={\varepsilon}^{-1}, \ \alpha_1=A_3, \ \alpha_3=A_0-A_2-{\varepsilon}^{-1}, \ \alpha_4=A_1
\end{gathered}\\
t=1+{\varepsilon}T, \ (x-1)(X-1)=1, \ (x-1)y+(X-1)Y=-A_2.
\end{gather*}
Notice that
$$
A_0+A_1+A_2+A_3=\alpha_0+\alpha_1+2\alpha_2+\alpha_3+\alpha_4=1
$$
and the change of variables from $(x,y)$ to $(X,Y)$ is symplectic.

As the fourth-order analogue of the above confluence process, we consider the following coupling confluence process from the Garnier system G(1,1,1,1,1) with symmetric Hamiltonians (see Theorem \ref{th:oka}). We take the following coupling confluence process $P_{VI} \rightarrow P_{V}$ for the coordinate system $(x,y)$ and $(z,w)$ of this system. Precisely speaking, for the Garnier system G(1,1,1,1,1), we make the change of parameters and variables
\begin{gather}
\begin{gathered}
\alpha_1=A_3+A_5-1, \quad \alpha_2=A_1, \quad \alpha_3=1-A_5,\\
\alpha_4=A_2, \quad \alpha_5=A_4+A_5-\frac{1}{\varepsilon}, \quad \alpha_6=1-A_3-A_5+\frac{1}{\varepsilon},
\end{gathered}\\
\begin{gathered}
T=\frac{t-1}{\varepsilon}, \quad  X=\frac{x}{x-1}, \quad Y=-(x-1)\{(x-1)y+\alpha_2\},\\
S=\frac{s-1}{\varepsilon}, \quad  Z=\frac{z}{z-1}, \quad W=-(z-1)\{(z-1)w+\alpha_4\}.
\end{gathered}
\end{gather}
from $\alpha_1,\alpha_2, \dots,\alpha_6,t,s,x,y,z,w$ to $A_1,\dots ,A_5,\varepsilon,T,S,X,Y,Z,W$. Then this system can also be written in the new variables $T,S,X,Y,Z,W$ and parameters $A_1,\dots,A_5,\varepsilon$ as a Hamiltonian system. This new system tends to the degenerate Garnier system G(1,1,1,2) with symmetric Hamiltonians (see Theorem \ref{symmetricDGV}) as $\varepsilon \rightarrow 0$ (see in Section 3).

We also show the confluence process among each system by taking the coupling confluence process of the Painlev\'e systems (see Figure 1).
\begin{figure}[ht]
\unitlength 0.1in
\begin{picture}( 35.9000, 30.8000)( 17.5000,-33.5000)
\put(17.5000,-8.6000){\makebox(0,0)[lb]{$P_{VI}$}}%
%
\special{pn 8}%
\special{pa 2290 780}%
\special{pa 2660 780}%
\special{fp}%
\special{sh 1}%
\special{pa 2660 780}%
\special{pa 2594 760}%
\special{pa 2608 780}%
\special{pa 2594 800}%
\special{pa 2660 780}%
\special{fp}%
\put(29.1000,-8.7000){\makebox(0,0)[lb]{$P_{V}$}}%
%
\special{pn 8}%
\special{pa 3440 720}%
\special{pa 3940 370}%
\special{fp}%
\special{sh 1}%
\special{pa 3940 370}%
\special{pa 3874 392}%
\special{pa 3896 402}%
\special{pa 3898 426}%
\special{pa 3940 370}%
\special{fp}%
%
\special{pn 8}%
\special{pa 3440 880}%
\special{pa 3930 1190}%
\special{fp}%
\special{sh 1}%
\special{pa 3930 1190}%
\special{pa 3884 1138}%
\special{pa 3886 1162}%
\special{pa 3864 1172}%
\special{pa 3930 1190}%
\special{fp}%
\put(42.1000,-4.4000){\makebox(0,0)[lb]{$P_{IV}$}}%
\put(42.3000,-12.7000){\makebox(0,0)[lb]{$P_{III}$}}%
%
\special{pn 8}%
\special{pa 4820 1190}%
\special{pa 5300 910}%
\special{fp}%
\special{sh 1}%
\special{pa 5300 910}%
\special{pa 5232 926}%
\special{pa 5254 938}%
\special{pa 5252 962}%
\special{pa 5300 910}%
\special{fp}%
\put(53.3000,-8.6000){\makebox(0,0)[lb]{$P_{II}$}}%
%
\special{pn 8}%
\special{pa 2430 1830}%
\special{pa 2430 2180}%
\special{fp}%
\special{sh 1}%
\special{pa 2430 2180}%
\special{pa 2450 2114}%
\special{pa 2430 2128}%
\special{pa 2410 2114}%
\special{pa 2430 2180}%
\special{fp}%
\put(25.9000,-21.0000){\makebox(0,0)[lb]{Coupling confluence process}}%
\put(17.5000,-28.3000){\makebox(0,0)[lb]{\eqref{uraS}}}%
%
\special{pn 8}%
\special{pa 2300 2740}%
\special{pa 2670 2740}%
\special{fp}%
\special{sh 1}%
\special{pa 2670 2740}%
\special{pa 2604 2720}%
\special{pa 2618 2740}%
\special{pa 2604 2760}%
\special{pa 2670 2740}%
\special{fp}%
\put(29.2000,-28.3000){\makebox(0,0)[lb]{\eqref{dVV}}}%
%
\special{pn 8}%
\special{pa 3440 2680}%
\special{pa 3940 2330}%
\special{fp}%
\special{sh 1}%
\special{pa 3940 2330}%
\special{pa 3874 2352}%
\special{pa 3896 2362}%
\special{pa 3898 2386}%
\special{pa 3940 2330}%
\special{fp}%
%
\special{pn 8}%
\special{pa 3440 2840}%
\special{pa 3930 3150}%
\special{fp}%
\special{sh 1}%
\special{pa 3930 3150}%
\special{pa 3884 3098}%
\special{pa 3886 3122}%
\special{pa 3864 3132}%
\special{pa 3930 3150}%
\special{fp}%
\put(42.3000,-24.0000){\makebox(0,0)[lb]{\eqref{SdeG3}}}%
\put(42.3000,-32.4000){\makebox(0,0)[lb]{\eqref{SdeGa}}}%
\put(53.4000,-28.3000){\makebox(0,0)[lb]{\eqref{SdeG4}}}%
\put(17.5000,-28.3000){\makebox(0,0)[lb]{\eqref{uraS}}}%
\put(20.4000,-15.0000){\makebox(0,0)[lb]{Confluence process of the Painlev\'e equations}}%
\put(20.0000,-35.2000){\makebox(0,0)[lb]{Confluence process of the Garnier system in two variables}}%
%
\special{pn 8}%
\special{pa 4750 360}%
\special{pa 5300 730}%
\special{fp}%
\special{sh 1}%
\special{pa 5300 730}%
\special{pa 5256 676}%
\special{pa 5256 700}%
\special{pa 5234 710}%
\special{pa 5300 730}%
\special{fp}%
%
\special{pn 8}%
\special{pa 4740 2330}%
\special{pa 5330 2700}%
\special{fp}%
\special{sh 1}%
\special{pa 5330 2700}%
\special{pa 5284 2648}%
\special{pa 5286 2672}%
\special{pa 5264 2682}%
\special{pa 5330 2700}%
\special{fp}%
\end{picture}%
\label{fig:degeneGarnier21}
\caption{}
\end{figure}

In this paper, we do not consider the coupling confluence process $P_{III} \rightarrow P_{II}$. These coupling confluence processes are new.

Moreover, we can make a generalization of these systems involving symmetry and holomorphy conditions without monodromy preserving deformations equations of second-order linear ordinary differential equations. The Hamiltonian system obtained by this process is generally rational functions of $t,s,x,y,z,w,q,p$, and is not symmetric. However, our Hamiltonians are all polynomials with respect to $x,y,z,w,q,p$ and symmetric. The system is explicitly given by
\begin{align}
\begin{split}
&dx=\frac{\partial H_1}{\partial y}dt+\frac{\partial H_2}{\partial y}ds+\frac{\partial H_3}{\partial y}du, \quad dy=-\frac{\partial H_1}{\partial x}dt-\frac{\partial H_2}{\partial x}ds-\frac{\partial H_3}{\partial x}du,\\
&dz=\frac{\partial H_1}{\partial w}dt+\frac{\partial H_2}{\partial w}ds+\frac{\partial H_3}{\partial w}du, \quad dw=-\frac{\partial H_1}{\partial z}dt-\frac{\partial H_2}{\partial z}ds-\frac{\partial H_3}{\partial z}du,\\
&dq=\frac{\partial H_1}{\partial p}dt+\frac{\partial H_2}{\partial p}ds+\frac{\partial H_3}{\partial p}du, \quad dp=-\frac{\partial H_1}{\partial q}dt-\frac{\partial H_2}{\partial q}ds-\frac{\partial H_3}{\partial q}du\\
\end{split}
\end{align}
with the symmetric polynomial Hamiltonians $H_i \in {\Bbb C}(t,s,u)[x,y,z,w,q,p] \ (i=1,2,3)$
\begin{align}\label{SdeGHHHH}
\begin{split}
&H_1:=H_{*}(x,y,t;\alpha_0,\ldots)+R(x,y,z,w,t,s;\alpha_0,\ldots)+R(x,y,q,p,t,u;\alpha_0,\ldots)\\
&H_2=\pi(H_1), \quad H_3=(\pi \circ \pi)(H_1),
\end{split}
\end{align}
where the transformation $\pi$ is explicitly given by
\begin{equation}
\pi:(x,y,z,w,q,p,t,s,u) \rightarrow (z,w,q,p,x,y,s,u,t), \quad (\pi)^3=1
\end{equation}
with some parameter's change, and the symbol $H_{*}(x,y,t;\alpha_0,\ldots)$ denotes one of the Painlev\'e Hamiltonians.

Moreover, we can make the symmetry and holomorphy conditions inductively. These Hamiltonian systems in three variables are new.

\section{Garnier system in two variables}
The Garnier system in two variables is equivalent to the Hamiltonian system
\begin{align}\label{1}
\begin{split}
&dx=\frac{\partial H_1}{\partial y}dt+\frac{\partial H_2}{\partial y}ds, \quad dy=-\frac{\partial H_1}{\partial x}dt-\frac{\partial H_2}{\partial x}ds,\\
&dz=\frac{\partial H_1}{\partial w}dt+\frac{\partial H_2}{\partial w}ds, \quad dw=-\frac{\partial H_1}{\partial z}dt-\frac{\partial H_2}{\partial z}ds
\end{split}
\end{align}
with the polynomial Hamiltonians
\begin{align}
\begin{split}
H_1 =&H_{VI}(x,y,t;1-2\alpha_1-\alpha_2-\alpha_3-\alpha_5,\alpha_2,\alpha_1,\alpha_5,\alpha_3)\\
&-\frac{\alpha_4 s}{t(t-s)}xy-\frac{\alpha_3(s-1)}{(t-1)(t-s)}zw+\frac{2(s-1)xyzw}{(t-1)(t-s)}\\
&-\frac{t(xy-\alpha_3)yz+s(zw-\alpha_4)xw}{t(t-s)}+\frac{\{2(xy+\alpha_1)+zw+\alpha_2\}xzw}{t(t-1)},
\end{split}\\
H_2=&\pi(H_1),
\end{align}
where the transformation $\pi$ is explicitly given by
\begin{align}
\begin{split}
\pi:(x,y,z,w,t,s;\alpha_1,\alpha_2,\alpha_3,\alpha_4,\alpha_5,\alpha_6) \rightarrow (z,w,x,y,s,t;\alpha_1,\alpha_2,\alpha_4,\alpha_3,\alpha_5,\alpha_6),
\end{split}
\end{align}
and the symbol $H_{VI}(x,y,t;\alpha_0,\alpha_1,\ldots,\alpha_4)$ denotes the sixth Painlev\'e Hamiltonian given by
\begin{align}
\begin{split}
&H_{VI}(x,y,t;\alpha_0,\alpha_1,\alpha_2,\alpha_3,\alpha_4)\\
&=\frac{1}{t(t-1)}[y^2(x-t)(x-1)x-\{(\alpha_0-1)(x-1)x+\alpha_3(x-t)x\\
&+\alpha_4(x-t)(x-1)\}y+\alpha_2(\alpha_1+\alpha_2)(x-t)] \quad (\alpha_0+\alpha_1+2\alpha_2+\alpha_3+\alpha_4=1).
\end{split}
\end{align}

\begin{theorem}
Let us consider a polynomial Hamiltonian system with Hamiltonian\\
$H_i \in {\Bbb C}(t,s)[x,y,z,w] \ (i=1,2)$. We assume that

$(A1)$ $deg(H_i)=5$ with respect to $x,y,z,w$.

$(A2)$ This system becomes again a polynomial Hamiltonian system in each coordinate system\\ $\{U_j,(x_j,y_j,z_j,w_j)\} \ (j=1,2,..,6)${\rm : \rm}
$$
{U_j}={{\Bbb C}^4} \ni {(x_j,y_j,z_j,w_j)} \ (j=1,2,..,6),
$$
via the following birational and symplectic transformations
\begin{align*}
&1) \ x_1=\frac{1}{x}, \ y_1=-x(xy+zw+\alpha_1), \ z_1=\frac{z}{x}, \ w_1=xw, \\
&2) \ x_2=\frac{1}{x}, \ y_2=-x(xy+zw+\alpha_1+\alpha_2), \ z_2=\frac{z}{x}, \ w_2=xw, \\
&3) \ x_3=-y(xy-\alpha_3),\ y_3=\frac{1}{y},\ z_3=z,\ w_3=w, \\
&4) \ x_4=x, \ y_4=y, \ z_4=-w(zw-\alpha_4), \ w_4=\frac{1}{w}, \\
&5) \ x_5=-((x+z-1)y-\alpha_5)y, \ y_5=\frac{1}{y}, \  z_5=z, \ w_5=w-y, \\
&6) \ x_6=-\left(\left(x+tz/s-t \right)y-\alpha_6 \right)y,\ y_6=\frac{1}{y}, \ z_6=z, \  w_6=w-ty/s.
\end{align*}
Then such a system coincides with the system $(1)$.
\end{theorem}

\begin{theorem} 
On each affine open set $(x_j,y_j,z_j,w_j) \in U_j \times B$ in Theorem 1.1, the  Hamiltonians $H_{j1}$ and $H_{j2}$ on $U_j \times B$ are expressed as a polynomial in $x_j,y_j,z_j,w_j$ and a rational function $t$ and $s$, and satisfies the following conditions{\rm: \rm}
\begin{align*}
\begin{split}
&dx \wedge dy +dz \wedge dw - dH_1 \wedge dt- dH_2 \wedge ds\\
&=dx_j \wedge dy_j +dz_j \wedge dw_j - dH_{j1} \wedge dt- dH_{j2} \wedge ds \quad (j=1,2,..,5),
\end{split}\\
\begin{split}
&dx \wedge dy +dz \wedge dw - d(H_1-(1-z/s)y) \wedge dt- d(H_2-(1-x/t)w) \wedge ds\\
&=dx_6 \wedge dy_6 +dz_6 \wedge dw_6 - dH_{61} \wedge dt- dH_{62} \wedge ds.
\end{split}
\end{align*}
\end{theorem}

\begin{theorem}\label{th:oka}
The birational and symplectic transformation
\begin{equation}
S:(x,y,z,w;t,s) \rightarrow (x+(zw+\alpha_1)/y,y,w/y,-zy,t,t/s),
\end{equation}
takes the system $(1)$ to the Hamiltonian system
\begin{align}\label{uraS}
\begin{split}
&dx=\frac{\partial H_1}{\partial y}dt+\frac{\partial H_2}{\partial y}ds, \quad dy=-\frac{\partial H_1}{\partial x}dt-\frac{\partial H_2}{\partial x}ds,\\
&dz=\frac{\partial H_1}{\partial w}dt+\frac{\partial H_2}{\partial w}ds, \quad dw=-\frac{\partial H_1}{\partial z}dt-\frac{\partial H_2}{\partial z}ds
\end{split}
\end{align}
with the polynomial Hamiltonians
\begin{align}\label{uraSH}
\begin{split}
H_1 &=H_{VI}(x,y,t;1-\alpha_1-\alpha_2-\alpha_3-\alpha_5,-\alpha_1-\alpha_2,\alpha_2,\alpha_1+\alpha_5,\alpha_1+\alpha_3)\\
&+\frac{\alpha_4 xy}{(t-s)}+\frac{\alpha_2(s-1)zw+2(s-1)xyzw}{(t-1)(t-s)}+\frac{\{ty+(xy+\alpha_2)x\}zw}{t(t-1)}\\
&-\frac{t(zw+\alpha_4)yz+s(xy+\alpha_2)xw}{t(t-s)},
\end{split}\\
H_2 &=\pi(H_1),
\end{align}
where the transformation $\pi$ is explicitly given by
\begin{equation}
\pi:(x,y,z,w;t,s,\alpha_1,\alpha_2, \dots ,\alpha_6) \rightarrow (z,w,x,y;s,t,\alpha_1,\alpha_4,\alpha_3,\alpha_2,\alpha_5,\alpha_6).
\end{equation}
\end{theorem}
This transformation can be regarded as a generalization of Okamoto-transformation of the sixth Painlev\'e system.

\begin{remark} {\rm
The system (1) is not invariant under $S$.
\rm}
\end{remark}

\begin{theorem}
Let us consider a polynomial Hamiltonian system with Hamiltonian\\
$H_i \in {\Bbb C}(t,s)[x,y,z,w] \ (i=1,2)$. We assume that

$(A1)$ $deg(H_i)=5$ with respect to $x,y,z,w$.

$(A2)$ This system becomes again a polynomial Hamiltonian system in each coordinate system $\{U_j,(x_j,y_j,z_j,w_j)\} \ (j=1,2,..,6)${\rm : \rm}
$$
{U_j}={{\Bbb C}^4} \ni {(x_j,y_j,z_j,w_j)} \ (j=1,2,..,6),
$$
via the following birational and symplectic transformations
\begin{align*}
&1)  x_1=\frac{1}{x}, \ y_1=-x(xy+\alpha_2), \ z_1=z, \ w_1=w,\\
&2)  x_2=1/x, \ y_2=-x(xy+zw-\alpha_1), \ z_2=\frac{z}{x}, \ w_2=xw,\\
&3)  x_3=x,\ y_3=y,\ z_3=1/z,\ w_3=-(zw+\alpha_4)z,\\
&4)  x_4=-(xy+zw-(\alpha_1+\alpha_3))y, \ y_4=1/y, \ z_4=zy, \ w_4=w/y,\\
&5)  x_5=-((x-1)y+(z-1)w-(\alpha_1+\alpha_5))y, \ y_5=1/y, \ z_5=(z-1)y, \ w_5=w/y,\\
&6)  x_6=-((x-t)y+(z-s)w-(\alpha_1+\alpha_6))y, \ y_6=1/y, \ z_6=(z-s)y, \ w_6=w/y.
\end{align*}
Then such a system coincides with the system \eqref{uraS}.
\end{theorem}

\begin{theorem} 
On each affine open set $(x_j,y_j,z_j,w_j) \in U_j \times B$ in Theorem \ref{th:oka}, the  Hamiltonians $H_{j1}$ and $H_{j2}$ on $U_j \times B$ are expressed as a polynomial in $x_j,y_j,z_j,w_j$ and a rational function in $t$ and $s$, and satisfies the following conditions{\rm: \rm}
\begin{align*}
\begin{split}
&dx \wedge dy +dz \wedge dw - dH_1 \wedge dt- dH_2 \wedge ds\\
&=dx_j \wedge dy_j +dz_j \wedge dw_j - dH_{j1} \wedge dt- dH_{j2} \wedge ds \quad (j=1,2,..,5),
\end{split}\\
\begin{split}
dx \wedge dy +dz \wedge dw - d(H_1-y) \wedge dt- d(H_2-w) \wedge ds\\
=dx_6 \wedge dy_6 +dz_6 \wedge dw_6 - dH_{61} \wedge dt- dH_{62} \wedge ds.
\end{split}
\end{align*}
\end{theorem}

\begin{theorem}\label{th:in3}
The system \eqref{uraS} is invariant under the following transformations\rm{: \rm} with the notation $(*)=(x,y,z,w,t,s;\alpha_1,\alpha_2,\ldots,\alpha_6),$
\begin{align*}
        u_1: (*) &\rightarrow (x+\alpha_2/y,y,z,w,t,s;\alpha_1+\alpha_1,-\alpha_2,\alpha_3,\alpha_4,\alpha_5,\alpha_6), \\
        u_2: (*) &\rightarrow (x,y,z+\alpha_4/w,w,t,s;\alpha_1+\alpha_4,\alpha_2,\alpha_3,-\alpha_4,\alpha_5,\alpha_6),\\
        u_3: (*) &\rightarrow (\frac{x(xy+zw-\alpha_1)}{(xy+zw-\alpha_1-\alpha_3)},\frac{y(xy+zw-\alpha_1-\alpha_3)}{(xy+zw-\alpha_1)},\frac{z(xy+zw-\alpha_1)}{(xy+zw-\alpha_1-\alpha_3)},\\
        &\frac{w(xy+zw-\alpha_1-\alpha_3)}{(xy+zw-\alpha_1)},t,s;\alpha_1+\alpha_3,\alpha_2,-\alpha_3,\alpha_4,\alpha_5,\alpha_6),\\
        \varphi_1: (*) &\rightarrow (1/x,-(xy+\alpha_2)x,1/z,-(zw+\alpha_4)z,1/t,1/s; \\
        &-\alpha_1-\alpha_2-\alpha_3-\alpha_4,\alpha_2,\alpha_3,\alpha_4,1-\alpha_6,1-\alpha_5),\\
        \varphi_2: (*) &\rightarrow (1-x,-y,1-z,-w,1-t,1-s;\alpha_1,\alpha_2,\alpha_5,\alpha_4,\alpha_3,\alpha_6),\\
        \varphi_3: (*) &\rightarrow \left(\frac{t-x}{t-1},-(t-1)y,\frac{s-z}{s-1},-(s-1)w,\frac{t}{t-1},\frac{s}{s-1};\alpha_1,\alpha_2,\alpha_6,\alpha_4,\alpha_5,\alpha_3 \right).
        \end{align*}
\end{theorem}

\section{A generalization of the system  \eqref{uraS} to three variables}

In this section, we present a generalization of the system \eqref{uraS} to three variables $t,s$ and $u$, which is equivalent to the polynomial Hamiltonian system
\begin{align}\label{SdeGH}
\begin{split}
&dx=\frac{\partial H_1}{\partial y}dt+\frac{\partial H_2}{\partial y}ds+\frac{\partial H_3}{\partial y}du, \quad dy=-\frac{\partial H_1}{\partial x}dt-\frac{\partial H_2}{\partial x}ds-\frac{\partial H_3}{\partial x}du,\\
&dz=\frac{\partial H_1}{\partial w}dt+\frac{\partial H_2}{\partial w}ds+\frac{\partial H_3}{\partial w}du, \quad dw=-\frac{\partial H_1}{\partial z}dt-\frac{\partial H_2}{\partial z}ds-\frac{\partial H_3}{\partial z}du,\\
&dq=\frac{\partial H_1}{\partial p}dt+\frac{\partial H_2}{\partial p}ds+\frac{\partial H_3}{\partial p}du, \quad dp=-\frac{\partial H_1}{\partial q}dt-\frac{\partial H_2}{\partial q}ds-\frac{\partial H_3}{\partial q}du\\
\end{split}
\end{align}
with the symmetric Hamiltonians $H_i \in {\Bbb C}(t,s,u)[x,y,z,w,q,p] \ (i=1,2,3)$
\begin{align}\label{SdeGHHHH55}
\begin{split}
H_1 &=H_{VI}(x,y,t;1-\alpha_1-\alpha_2-\alpha_3-\alpha_5,-\alpha_1-\alpha_2,\alpha_2,\alpha_1+\alpha_5,\alpha_1+\alpha_3)\\
&+R(x,y,z,w,t,s;\alpha_2,\alpha_4)+R(x,y,q,p,t,u;\alpha_2,\alpha_7),\\
H_2 &=\pi(H_1), \quad H_3=(\pi \circ \pi)(H_1) \quad (2\alpha_1+\alpha_2+ \dots +\alpha_7=1),
\end{split}
\end{align}
where the transformation $\pi$ is explicitly given by
\begin{equation}
\pi:(*) \rightarrow (z,w,q,p,x,y,s,u,t;\alpha_1,\alpha_4,\alpha_3,\alpha_7,\alpha_5,\alpha_6,\alpha_2).
\end{equation}
Here the symbol $(*)$ denotes $(*):=(x,y,z,w,q,p,t,s,u;\alpha_1,\alpha_2,\alpha_3,\alpha_4,\alpha_5,\alpha_6,\alpha_7)$, and the symbol\\
 $R(q_l,p_l,q_m,p_m,t_l,t_m;\alpha,\beta)$ is explicitly given by
\begin{align*}
&R(q_l,p_l,q_m,p_m,t_l,t_m;\alpha,\beta)\\
&=\frac{\beta q_lp_l}{t_l-t_m}+\frac{\alpha (t_m-1)q_mp_m}{(t_l-1)(t_l-t_m)}+\frac{\{t_l p_l+(q_lp_l+\alpha)q_l\}q_mp_m}{t_l(t_l-1)}\\
&-\frac{t_l(q_mp_m+\beta)p_lq_m+t_m(q_lp_l+\alpha)q_lp_m}{t_l(t_l-t_m)}.
\end{align*}
For the Hamiltonian system \eqref{1}, it is difficult to make a generalization to three variables by a similar way.

\begin{theorem}
Let us consider a polynomial Hamiltonian system with Hamiltonian\\
$H_i \in {\Bbb C}(t,s,u)[x,y,z,w,q,p] \ (i=1,2,3)$. We assume that

$(A1)$ $deg(H_i)=5$ with respect to $x,y,z,w,q,p$.

$(A2)$ This system becomes again a polynomial Hamiltonian system in each coordinate system $\{U_j,(x_j,y_j,z_j,w_j,q_j,p_j)\} \ (j=1,2, \dots,7)${\rm : \rm}
$$
{U_j}={{\Bbb C}^6} \ni {(x_j,y_j,z_j,w_j,q_j,p_j)} \ (j=1,2, \dots,7),
$$
via the following birational and symplectic transformations
\begin{align*}
&1)  x_1=1/x, \ y_1=-x(xy+zw+qp-\alpha_1), \ z_1=\frac{z}{x}, \ w_1=xw, \ q_1=\frac{q}{x}, \ p_1=xp,\\
&2)  x_2=\frac{1}{x}, \ y_2=-x(xy+\alpha_2), \ z_2=z, \ w_2=w, \ q_2=q, \ p_2=p,\\
&3)  x_3=-(xy+zw+qp-(\alpha_1+\alpha_3))y, \ y_3=1/y, \ z_3=zy, \ w_3=w/y,\\
&q_3=qy, \ w_3=p/y,\\
&4)  x_4=x,\ y_4=y,\ z_4=1/z,\ w_4=-(zw+\alpha_4)z, \ q_4=q, \ p_4=p,\\
&5)  x_5=-((x-1)y+(z-1)w+(q-1)p-(\alpha_1+\alpha_5))y, \ y_5=1/y,\\
&z_5=(z-1)y, \ w_5=w/y, \ q_5=(q-1)y, \ p_5=p/y,\\
&6)  x_6=-((x-t)y+(z-s)w+(q-u)p-(\alpha_1+\alpha_6))y, \ y_6=1/y,\\
&z_6=(z-s)y, \ w_6=w/y, \ q_6=(q-u)y, \ p_6=p/y,\\
&7)  x_7=x, \ y_7=y, \ z_7=z, \ w_7=w, \ q_7=1/q, \ p_7=-(qp+\alpha_7)q.
\end{align*}
Then such a system coincides with the system \eqref{SdeGHHHH55}.
\end{theorem}

\begin{theorem}
The system \eqref{SdeGHHHH55} is invariant under the following transformations\rm{:\rm} with the notation $(*)=(x,y,z,w,q,p,t,s,u;\alpha_1,\alpha_2,\ldots,\alpha_7),$
\begin{align*}
        u_2: (*) &\rightarrow (x+\alpha_2/y,y,z,w,q,p,t,s,u;\alpha_1+\alpha_1,-\alpha_2,\alpha_3,\alpha_4,\alpha_5,\alpha_6,\alpha_7), \\
        u_4: (*) &\rightarrow (x,y,z+\alpha_4/w,w,q,p,t,s,u;\alpha_1+\alpha_4,\alpha_2,\alpha_3,-\alpha_4,\alpha_5,\alpha_6,\alpha_7),\\
        u_7: (*) &\rightarrow (x,y,z,w,q,p+\alpha_7/p,t,s,u;\alpha_1+\alpha_7,\alpha_2,\alpha_3,\alpha_4,\alpha_5,\alpha_6,-\alpha_7),\\
        u_3: (*) &\rightarrow (\frac{x(xy+zw+qp-\alpha_1)}{(xy+zw+qp-\alpha_1-\alpha_3)},\frac{y(xy+zw+qp-\alpha_1-\alpha_3)}{(xy+zw+qp-\alpha_1)},\\
        &\frac{z(xy+zw+qp-\alpha_1)}{(xy+zw+qp-\alpha_1-\alpha_3)},\frac{w(xy+zw+qp-\alpha_1-\alpha_3)}{(xy+zw+qp-\alpha_1)},\\
        &\frac{q(xy+zw+qp-\alpha_1)}{(xy+zw+qp-\alpha_1-\alpha_3)},\frac{p(xy+zw+qp-\alpha_1-\alpha_3)}{(xy+zw+qp-\alpha_1)},\\
        &t,s,u;\alpha_1+\alpha_3,\alpha_2,-\alpha_3,\alpha_4,\alpha_5,\alpha_6,\alpha_7),\\
        \varphi_1: (*) &\rightarrow (1/x,-(xy+\alpha_2)x,1/z,-(zw+\alpha_4)z,1/q,-(qp+\alpha_7)q,1/t,1/s,1/u; \\
        &-\alpha_1-\alpha_2-\alpha_3-\alpha_4-\alpha_7,\alpha_2,\alpha_3,\alpha_4,1-\alpha_6,1-\alpha_5,\alpha_7),\\
        \varphi_2: (*) &\rightarrow (1-x,-y,1-z,-w,1-q,-p,1-t,1-s,1-u;\\
        &\alpha_1,\alpha_2,\alpha_5,\alpha_4,\alpha_3,\alpha_6,\alpha_7),\\
        \varphi_3: (*) &\rightarrow (\frac{t-x}{t-1},-(t-1)y,\frac{s-z}{s-1},-(s-1)w,\frac{u-q}{u-1},-(u-1)p,\\
        &\frac{t}{t-1},\frac{s}{s-1},\frac{u}{u-1};\alpha_1,\alpha_2,\alpha_6,\alpha_4,\alpha_5,\alpha_3,\alpha_7).
        \end{align*}
\end{theorem}

\section{Degenerate Garnier system $G(1,1,1,2)$ in two variables}

The degenerate Garnier system $G(1,1,1,2)$ in two variables $t,s$ is equivalent to the Hamiltonian system
\begin{align}\label{dV}
\begin{split}
&dx=\frac{\partial H_1}{\partial y}dt+\frac{\partial H_2}{\partial y}ds, \quad dy=-\frac{\partial H_1}{\partial x}dt-\frac{\partial H_2}{\partial x}ds,\\
&dz=\frac{\partial H_1}{\partial w}dt+\frac{\partial H_2}{\partial w}ds, \quad dw=-\frac{\partial H_1}{\partial z}dt-\frac{\partial H_2}{\partial z}ds
\end{split}
\end{align}
with the polynomial Hamiltonians $K_i \in {\Bbb C}(t,s)[x,y,z,w] \ (i=1,2)$ (see \cite{K2})
\begin{align}\label{K1}
\begin{split}
t^2K_1 =&x^2(x-t)y^2+2x^2zyw+xz(z-s)w^2\\
&-\{(\alpha_0+\alpha_2-1)x^2+\alpha_1x(x-t)+\eta(x-t)+\eta tz\}y\\
&-\{(\alpha_0+\alpha_1-1)xz+\alpha_2x(z-s)+\eta(s-1)z\}w+\nu(\nu+\alpha_3)x,\\
s(s-1)K_2 =&x^2zy^2+2xz(z-s)yw\\
&+\{z(z-1)(z-s)+\frac{s(s-1)}{t}xz\}w^2\\
&-\{(\alpha_0+\alpha_1-1)xz+\alpha_2x(z-s)-\eta(s-1)z\}y\\
&-\{(\alpha_0-1)z(z-1)+\alpha_1z(z-s)+\alpha_2(z-1)(z-s)\\
&+\frac{s(s-1)}{t}(\alpha_2x+\eta z)\}w+\nu(\nu+\alpha_3)z \ \left(\nu=-\frac{1}{2}(\alpha_0+\alpha_1+\alpha_2-1+\alpha_3) \right).
\end{split}
\end{align}

\begin{theorem}
The system \eqref{dV} is invariant under the following transformations\rm{:\rm} with the notation $(*)=(x,y,z,w,\eta,t,s;\alpha_0,\alpha_1,\ldots,\alpha_3,\nu),$
\begin{align*}
        s_0: (*) &\rightarrow \left(x,y-\frac{s\alpha_0}{sx+tz-ts},z,w-\frac{t\alpha_0}{sx+tz-ts},\eta,t,s;-\alpha_0,\alpha_1,\alpha_2,-\alpha_3,\nu \right),\\
        s_1: (*) &\rightarrow \left(x,y-\frac{\alpha_1}{x}+\frac{\eta(z-1)}{x^2},z,w-\frac{\eta}{x},-\eta,t,s;\alpha_0,-\alpha_1,\alpha_2,-\alpha_3,\nu \right),\\
        s_2: (*) &\rightarrow \left(x,y,z,w-\frac{\alpha_2}{z},\eta,t,s;\alpha_0,\alpha_1,-\alpha_2,\alpha_3,\nu+\alpha_2 \right),\\
        s_3: (*) &\rightarrow (x,y,z,w,\eta,t,s;\alpha_0,\alpha_1,\alpha_2,-\alpha_3,\nu+\alpha_3).
\end{align*}
\end{theorem}

\begin{theorem}
Let us consider a polynomial Hamiltonian system with Hamiltonian\\
$H_i \in {\Bbb C}(t,s)[x,y,z,w] \ (i=1,2)$. We assume that

$(A1)$ $deg(H_i)=5$ with respect to $x,y,z,w$.

$(A2)$ This system becomes again a polynomial Hamiltonian system in each coordinate system $\{U_j,(x_j,y_j,z_j,w_j)\} \ (j=1,2,..,5)${\rm : \rm}
$$
{U_j}={{\Bbb C}^4} \ni {(x_j,y_j,z_j,w_j)} \ (j=1,2,..,5),
$$
via the following birational and symplectic transformations
\begin{align*}
&1)x_1=\frac{1}{x}, \quad y_1=-(xy+zw+\nu+\alpha_3)x, \quad z_1=\frac{z}{x}, \quad w_1=xw,\\
&2)x_2=\frac{1}{x}, \quad y_2=-(xy+zw+\nu)x, \quad z_2=\frac{z}{x}, \quad w_2=xw,\\
&3)x_3=x,\quad y_3=y-\frac{\alpha_1}{x}+\frac{\eta(z-1)}{x^2},\quad z_3=z,\quad w_3=w-\frac{\eta}{x},\\
&4)x_4=x, \quad y_4=y, \quad z_4=-(zw-\alpha_2)w, \quad w_4=\frac{1}{w},\\
&5)x_5=-\left(\left(x+\frac{t}{s}z-t \right)y-\alpha_0 \right)y, \quad y_5=\frac{1}{y}, \quad  z_5=z, \quad w_5=w-\frac{t}{s}y.
\end{align*}
Then such a system coincides with the system \eqref{dV}.
\end{theorem}

\section{Symmetric Hamiltonian of the system \eqref{dV}}
In this section, we make {\it symmetric Hamiltonian} for the degenerate Garnier system $G(1,1,1,2)$ by taking suitable birational and symplectic transformations. By making this symmetric Hamiltonian, we can make a generalization of this system involving symmetry and holomorphy in the next section. We also show the confluence process from the system \eqref{uraS} to this system by taking the coupling confluence process $P_{VI} \rightarrow P_{V}$ for each coordinate system $(x,y)$ and $(z,w)$ of the system \eqref{uraS}, respectively.

\begin{theorem}\label{symmetricDGV}
The birational and symplectic transformations with parameter's change\rm{:\rm}
\begin{align}
\begin{split}
&X:=x(xy+zw+\nu), \quad Y:=\frac{1}{x}, \quad Z:=xw, \quad W:=-\frac{z}{x}, \quad T:=-\frac{1}{t}, \quad S:=-\frac{s}{t},\\
&\eta=1, \quad (\alpha_0,\alpha_1,\alpha_2,\alpha_3,\nu) \rightarrow (\alpha_5,\alpha_4-\alpha_3,\alpha_2,\alpha_1,\alpha_3)
\end{split}
\end{align}
takes the system \eqref{dV} to the polynomial Hamiltonian system
\begin{align}\label{dVV}
\begin{split}
&dx=\frac{\partial H_1}{\partial y}dt+\frac{\partial H_2}{\partial y}ds, \quad dy=-\frac{\partial H_1}{\partial x}dt-\frac{\partial H_2}{\partial x}ds,\\
&dz=\frac{\partial H_1}{\partial w}dt+\frac{\partial H_2}{\partial w}ds, \quad dw=-\frac{\partial H_1}{\partial z}dt-\frac{\partial H_2}{\partial z}ds
\end{split}
\end{align}
with the symmetric Hamiltonians $H_i \in {\Bbb C}(t,s)[x,y,z,w] \ (i=1,2)$
\begin{align}\label{SdeG}
\begin{split}
H_1 &=H_{V}(x,y,t;\alpha_3,\alpha_1,\alpha_4)-\frac{\alpha_2sxy}{t(s-t)}-\frac{\alpha_1zw}{s-t}\\
&+\frac{tx^2yw+tyzw+\alpha_1txw+syz^2w-syzw+\alpha_2syz-2txyzw}{t(s-t)},\\
H_2 &=\pi(H_1) \quad (\alpha_1+\alpha_2+ \dots +\alpha_5=1),
\end{split}
\end{align}
\end{theorem}
Here, for notational convenience, we have renamed $X,Y,Z,W,T,S$ to $x,y,z,w,t,s$ (which are not the same as the previous $x,y,z,w,t,s$). The transformation $\pi$ is explicitly given by
\begin{equation}
\pi:(x,y,z,w,t,s;\alpha_1,\alpha_2,\alpha_3,\alpha_4,\alpha_5)
 \rightarrow (z,w,x,y,s,t;\alpha_2,\alpha_1,\alpha_3,\alpha_4,\alpha_5).
\end{equation}
We note that the Hamiltonian $H_1$ involves the Painlev\'e V Hamiltonian $H_{V}$ given by
\begin{align}
\begin{split}
&H_{V}(x,y,t;\alpha_1,\alpha_2,\alpha_3)=\frac{x(x-1)y(y+t)+\alpha_2tx-\alpha_3xy-\alpha_1y(x-1)}{t}.
\end{split}
\end{align}

\begin{theorem}
Let us consider a polynomial Hamiltonian system with Hamiltonian\\
$H_i \in {\Bbb C}(t,s)[x,y,z,w] \ (i=1,2)$. We assume that

$(A1)$ $deg(H_i)=5$ with respect to $x,y,z,w$.

$(A2)$ This system becomes again a polynomial Hamiltonian system in each coordinate system $\{U_j,(x_j,y_j,z_j,w_j)\} \ (j=1,2,..,5)${\rm : \rm}
$$
{U_j}={{\Bbb C}^4} \ni {(x_j,y_j,z_j,w_j)} \ (j=1,2,..,5),
$$
via the following birational and symplectic transformations
\begin{align*}
&1)x_1=\frac{1}{x}, \quad y_1=-(yx+\alpha_1)x, \quad z_1=z, \quad w_1=w,\\
&2)x_2=x, \quad y_2=y, \quad z_2=\frac{1}{z}, \quad w_2=-(zw+\alpha_2)z,\\
&3)x_3=-(yx+wz-\alpha_3)y,\quad y_3=\frac{1}{y},\quad z_3=zy,\quad w_3=\frac{w}{y},\\
&4)x_4=-((x-1)y+(z-1)w-\alpha_4)y, \quad y_4=\frac{1}{y}, \quad z_4=(z-1)y, \quad w_4=\frac{w}{y},\\
&5)x_5=\frac{1}{x}, \quad y_5=-((y+tw/s+t)x+\alpha_5)x, \quad  z_5=z-tx/s, \quad w_5=w.
\end{align*}
Then such a system coincides with the system \eqref{dVV}.
\end{theorem}

\begin{theorem}
The system \eqref{dVV} is invariant under the following transformations\rm{:\rm} with the notation $(*)=(x,y,z,w,t,s;\alpha_1,\alpha_2,\ldots,\alpha_5),$
\begin{align*}
        s_1: (*) &\rightarrow \left(x+\frac{\alpha_1}{y},y,z,w,t,s;-\alpha_1,\alpha_2,\alpha_3+\alpha_1,\alpha_4+\alpha_1,\alpha_5 \right),\\
        s_2: (*) &\rightarrow \left(x,y,z+\frac{\alpha_2}{w},w,t,s;\alpha_1,-\alpha_2,\alpha_3+\alpha_2,\alpha_4+\alpha_2,\alpha_5 \right),\\
        s_5: (*) &\rightarrow \left(x+\frac{\alpha_5}{y+tw/s+t},y,z+\frac{t\alpha_5}{s(y+tw/s+t)},w,t,s;\alpha_1,\alpha_2,\alpha_3+\alpha_5,\alpha_4+\alpha_5,-\alpha_5 \right),\\
        \pi_1: (*) &\rightarrow (z,w,x,y,s,t;\alpha_2,\alpha_1,\alpha_3,\alpha_4,\alpha_5),\\
        \pi_2: (*) &\rightarrow (1-x,-y,1-z,-w,-t,-s;\alpha_1,\alpha_2,\alpha_4,\alpha_3,\alpha_5),\\
        \pi_3: (*) &\rightarrow (-\frac{xy+zw-\alpha_3}{tx},\frac{tx(xy+\alpha_1)}{xy+zw-\alpha_3},-\frac{xy+zw-\alpha_3}{sz},\frac{sz(zw+\alpha_2)}{xy+zw-\alpha_3},-t,-s;\\
        &\alpha_1,\alpha_2,\alpha_4+\alpha_5-1,\alpha_3+\alpha_5,1-\alpha_5),\\
        \pi_4: (*) &\rightarrow (-\frac{y(x-1)+w(z-1)-\alpha_4}{t(x-1)},\frac{t(x-1)((x-1)y+\alpha_1)}{y(x-1)+w(z-1)-\alpha_4},\\
        &-\frac{y(x-1)+w(z-1)-\alpha_4}{s(z-1)},\frac{s(z-1)((z-1)w+\alpha_2)}{y(x-1)+w(z-1)-\alpha_4},t,s;\\
        &\alpha_1,\alpha_2,\alpha_3+\alpha_5-1,\alpha_4+\alpha_5,1-\alpha_5).
\end{align*}
\end{theorem}

\begin{proposition}
The transformation $\pi_3$ can be obtained by composing the following transformations:

{\bf Step 1:} We transform the system \eqref{dVV} by the following birational and symplectic transformation
$$
g_1:(x,y,z,w) \rightarrow (x,y+(zw-\alpha_3)/x,z/x,xw).
$$

{\bf Step 2:} We then transform the system obtained by Step 1 by the following birational and symplectic transformation
$$
g_2:(x,y,z,w) \rightarrow (x-(zw-\alpha_1-\alpha_3)/y,y,z/y,wy).
$$

{\bf Step 3:} We then transform the system obtained by Step 2 by the following birational and symplectic transformation
$$
g_3:(*) \rightarrow (y,-x,1/z,-(zw+\alpha_2)z).
$$

{\bf Step 4:} We then transform the system obtained by Step 3 by the following birational and symplectic transformation
$$
g_4:(*) \rightarrow (-x/t,-ty,-z/s,-sw,-t,-s).
$$
\end{proposition}

As is well-known, the degeneration from $P_{VI}$ to $P_{V}$ (see \cite{T2}) is given by
\begin{gather*}
\begin{gathered}
\alpha_0={\varepsilon}^{-1}, \ \alpha_1=A_3, \ \alpha_3=A_0-A_2-{\varepsilon}^{-1}, \ \alpha_4=A_1
\end{gathered}\\
t=1+{\varepsilon}T, \ (x-1)(X-1)=1, \ (x-1)y+(X-1)Y=-A_2.
\end{gather*}
Notice that
$$
A_0+A_1+A_2+A_3=\alpha_0+\alpha_1+2\alpha_2+\alpha_3+\alpha_4=1
$$
and the change of variables from $(x,y)$ to $(X,Y)$ is symplectic.

As the fourth-order analogue of the above confluence process, we consider the following coupling confluence process from the system \eqref{uraS}. We take the following coupling confluence process $P_{VI} \rightarrow P_{V}$ for each coordinate system $(x,y)$ and $(z,w)$ of the system \eqref{uraS}.
\begin{theorem}
For the system \eqref{uraS}, we make the change of parameters and variables
\begin{gather}
\begin{gathered}\label{11}
\alpha_1=A_3+A_5-1, \quad \alpha_2=A_1, \quad \alpha_3=1-A_5,\\
\alpha_4=A_2, \quad \alpha_5=A_4+A_5-\frac{1}{\varepsilon}, \quad \alpha_6=1-A_3-A_5+\frac{1}{\varepsilon},
\end{gathered}\\
\begin{gathered}\label{12}
T=\frac{t-1}{\varepsilon}, \quad S=\frac{s-1}{\varepsilon}, \quad X=\frac{x}{x-1}, \quad Z=\frac{z}{z-1},\\
Y=-(x-1)\{(x-1)y+\alpha_2\}, \quad W=-(z-1)\{(z-1)w+\alpha_4\}
\end{gathered}
\end{gather}
from $\alpha_1,\alpha_2, \dots,\alpha_6,t,s,x,y,z,w$ to $A_1,\dots ,A_5,\varepsilon,T,S,X,Y,Z,W$. Then this system can also be written in the new variables $T,S,X,Y,Z,W$ and parameters $A_1,\dots,A_5,\varepsilon$ as a Hamiltonian system. This new system tends to the system \eqref{dVV} as $\varepsilon \rightarrow 0$.
\end{theorem}

\section{A generalization of the system \eqref{dVV} to three variables}
In this section, we present a generalization of the system \eqref{dVV} to three variables $t,s$ and $u$, which is equivalent to the polynomial Hamiltonian system
\begin{align}\label{deGHS}
\begin{split}
&dx=\frac{\partial H_1}{\partial y}dt+\frac{\partial H_2}{\partial y}ds+\frac{\partial H_3}{\partial y}du, \quad dy=-\frac{\partial H_1}{\partial x}dt-\frac{\partial H_2}{\partial x}ds-\frac{\partial H_3}{\partial x}du,\\
&dz=\frac{\partial H_1}{\partial w}dt+\frac{\partial H_2}{\partial w}ds+\frac{\partial H_3}{\partial w}du, \quad dw=-\frac{\partial H_1}{\partial z}dt-\frac{\partial H_2}{\partial z}ds-\frac{\partial H_3}{\partial z}du,\\
&dq=\frac{\partial H_1}{\partial p}dt+\frac{\partial H_2}{\partial p}ds+\frac{\partial H_3}{\partial p}du, \quad dp=-\frac{\partial H_1}{\partial q}dt-\frac{\partial H_2}{\partial q}ds-\frac{\partial H_3}{\partial q}du\\
\end{split}
\end{align}
with the symmetric Hamiltonians $H_i \in {\Bbb C}(t,s,u)[x,y,z,w,q,p] \ (i=1,2,3)$
\begin{align}\label{deGHSH}
\begin{split}
H_1 &=H_{V}(x,y,t;\alpha_4,\alpha_1,\alpha_5)\\
&+R(x,y,z,w,t,s;\alpha_1,\alpha_2)+R(x,y,q,p,t,u;\alpha_1,\alpha_3),\\
H_2 &=\pi(H_1), \quad H_3=(\pi \circ \pi)(H_1) \quad (\alpha_1+\alpha_2+ \dots +\alpha_6=1),
\end{split}
\end{align}
where the transformation $\pi$ is explicitly given by
\begin{equation}
\pi:(*) \rightarrow (z,w,q,p,x,y,s,u,t;\alpha_2,\alpha_3,\alpha_1,\alpha_4,\alpha_5,\alpha_6).
\end{equation}
Here the symbol $(*)$ denotes $(*):=(x,y,z,w,q,p,t,s,u;\alpha_1,\ldots,\alpha_6)$, and the symbol \\$R(q_l,p_l,q_m,p_m,t_l,t_m;\alpha,\beta)$ is explicitly given by
\begin{align*}
&R(q_l,p_l,q_m,p_m,t_l,t_m;\alpha,\beta)\\
&=\frac{\beta t_m q_lp_l}{t_l(t_l-t_m)}+\frac{\alpha q_mp_m}{t_l-t_m}\\
&-\frac{t_l{q_l}^2p_lp_m+t_lp_lq_mp_m+\alpha t_lq_lp_m+t_mp_l{q_m}^2p_m-t_mp_lq_mp_m+\beta t_m p_lq_m-2t_l q_lp_lq_mp_m}{t_l(t_l-t_m)}.
\end{align*}

\begin{theorem}
Let us consider a polynomial Hamiltonian system with Hamiltonian\\
$H_i \in {\Bbb C}(t,s,u)[x,y,z,w,q,p] \ (i=1,2,3)$. We assume that

$(A1)$ $deg(H_i)=5$ with respect to $x,y,z,w,q,p$.

$(A2)$ This system becomes again a polynomial Hamiltonian system in each coordinate system\\$\{U_j,(x_j,y_j,z_j,w_j,q_j,p_j)\} \ (j=1,2,..,6)${\rm : \rm}
$$
{U_j}={{\Bbb C}^6} \ni {(x_j,y_j,z_j,w_j,q_j,p_j)} \ (j=1,2,..,6),
$$
via the following birational and symplectic transformations
\begin{align*}
&1)x_1=\frac{1}{x}, \quad y_1=-(yx+\alpha_1)x, \quad z_1=z, \quad w_1=w, \quad q_1=q, \quad p_1=p,\\
&2)x_2=x, \quad y_2=y, \quad z_2=\frac{1}{z}, \quad w_2=-(zw+\alpha_2)z, \quad q_2=q, \quad p_2=p,\\
&3)x_3=x, \quad y_3=y, \quad z_3=z, \quad w_3=w, \quad q_3=\frac{1}{q}, \quad p_3=-(qp+\alpha_3)q,\\
&4)x_4=-(yx+wz+pq-\alpha_4)y,\quad y_4=\frac{1}{y},\quad z_4=zy,\\
&w_4=\frac{w}{y}, \quad q_4=qy, \quad p_4=\frac{p}{y},\\
&5)x_5=-((x-1)y+(z-1)w+(q-1)p-\alpha_5)y, \quad y_5=\frac{1}{y}, \quad z_5=(z-1)y,\\
&w_5=\frac{w}{y}, \quad q_5=(q-1)y, \quad p_5=\frac{p}{y},\\
&6)x_6=\frac{1}{x}, \quad y_6=-((y+tw/s+tp/u+t)x+\alpha_6)x, \quad  z_6=z-tx/s,\\
&w_6=w, \quad q_6=q-tx/u, \quad p_6=p.
\end{align*}
Then such a system coincides with the system \eqref{deGHS}.
\end{theorem}

\begin{theorem}
The system \eqref{deGHS} is invariant under the following transformations\rm{:\rm} with the notation $(*)=(x,y,z,w,q,p,t,s,u;\alpha_1,\alpha_2,\ldots,\alpha_6),$
\begin{align*}
        s_1: (*) &\rightarrow \left(x+\frac{\alpha_1}{y},y,z,w,q,p,t,s,u;-\alpha_1,\alpha_2,\alpha_3,\alpha_4+\alpha_1,\alpha_5+\alpha_1,\alpha_6 \right),\\
        s_2: (*) &\rightarrow \left(x,y,z+\frac{\alpha_2}{w},w,t,s,u;\alpha_1,-\alpha_2,\alpha_3,\alpha_4+\alpha_2,\alpha_5+\alpha_2,\alpha_6 \right),\\
        s_3: (*) &\rightarrow \left(x,y,z,w,q+\frac{\alpha_3}{p},p,t,s,u;\alpha_1,\alpha_2,-\alpha_3,\alpha_4+\alpha_3,\alpha_5+\alpha_3,\alpha_6 \right),\\
        s_6: (*) &\rightarrow (x+\frac{\alpha_6 su}{suy+tuw+tsp+tsu},y,z+\frac{\alpha_6 tu}{suy+tuw+tsp+tsu},w,\\
        &q+\frac{\alpha_6 ts}{suy+tuw+tsp+tsu},p,t,s,u;\alpha_1,\alpha_2,\alpha_3,\alpha_4+\alpha_6,\alpha_5+\alpha_6,-\alpha_6),\\
        \pi_1: (*) &\rightarrow (z,w,q,p,x,y,s,u,t;\alpha_3,\alpha_2,\alpha_1,\alpha_4,\alpha_5,\alpha_6),\\
        \pi_2: (*) &\rightarrow (1-x,-y,1-z,-w,1-q,-p,-t,-s,-u;\alpha_1,\alpha_2,\alpha_3,\alpha_5,\alpha_4,\alpha_6),\\
        \pi_3: (*) &\rightarrow (-\frac{xy+zw+qp-\alpha_4}{tx},\frac{tx(xy+\alpha_1)}{xy+zw+qp-\alpha_4},-\frac{xy+zw+qp-\alpha_4}{sz},\\
        &\frac{sz(zw+\alpha_2)}{xy+zw+qp-\alpha_4},-\frac{xy+zw+qp-\alpha_4}{uq},\frac{uq(qp+\alpha_3)}{xy+zw+qp-\alpha_4},-t,-s,-u;\\
        &\alpha_1,\alpha_2,\alpha_3,\alpha_5+\alpha_6-1,\alpha_4+\alpha_6,1-\alpha_6),\\
        \pi_4: (*) &\rightarrow (\frac{(1-x)y+(1-z)w+(1-q)p+\alpha_5}{t(x-1)},\\
        &\frac{t(x-1)((x-1)y+\alpha_1)}{(x-1)y+(z-1)w+(q-1)p-\alpha_5},\frac{(1-x)y+(1-z)w+(1-q)p+\alpha_5}{s(z-1)},\\
        &\frac{s(z-1)((z-1)w+\alpha_2)}{(x-1)y+(z-1)w+(q-1)p-\alpha_5},\frac{(1-x)y+(1-z)w+(1-q)p+\alpha_5}{u(q-1)},\\
        &\frac{u(q-1)((q-1)p+\alpha_3)}{(x-1)y+(z-1)w+(q-1)p-\alpha_5},t,s,u;\alpha_1,\alpha_2,\alpha_3,\alpha_4+\alpha_6-1,\alpha_5+\alpha_6,1-\alpha_6).
\end{align*}
\end{theorem}

\section{Degeneration from the system \eqref{dVV}}
As the fourth-order analogue of the confluence process from $P_{V}$ to $P_{IV}$ (see \cite{T2}), we consider the following coupling confluence process from the system \eqref{dVV}. We take the following coupling confluence process $P_{V} \rightarrow P_{IV}$ for each coordinate system $(x,y)$ and $(z,w)$ of the system \eqref{dVV}.
\begin{theorem}
For the system \eqref{dVV}, we make the change of parameters and variables
\begin{gather}
\begin{gathered}
\alpha_1=A_1, \quad \alpha_2=A_2, \quad \alpha_3=A_3, \quad \alpha_4=-\frac{1}{2{\varepsilon}^2}, \quad \alpha_5=1-A_1-A_2-A_3+\frac{1}{2{\varepsilon}^2}
\end{gathered}\\
\begin{gathered}
t=\frac{1+2\varepsilon T}{2{\varepsilon}^2}, \quad s=\frac{1+2\varepsilon S}{2{\varepsilon}^2}, \quad x=\frac{\varepsilon X}{\varepsilon X-1}, \quad z=\frac{\varepsilon Z}{\varepsilon Z-1},\\
y=-\frac{(\varepsilon X-1)\{(\varepsilon X-1)Y+\varepsilon A_2\}}{\varepsilon}, \quad w=-\frac{(\varepsilon Z-1)\{(\varepsilon Z-1)W+\varepsilon A_2\}}{\varepsilon}
\end{gathered}
\end{gather}
from $\alpha_1,\alpha_2, \dots,\alpha_5,t,x,y,z,w$ to $A_1,\dots ,A_4,\varepsilon,T,X,Y,Z,W$. Then this system can also be written in the new variables $T,X,Y,Z,W$ and parameters $A_1,A_2,A_3,A_4,\varepsilon$ as a Hamiltonian system. This new system tends to the polynomial Hamiltonian system 
\begin{align}\label{SdeG3}
\begin{split}
&dx=\frac{\partial H_1}{\partial y}dt+\frac{\partial H_2}{\partial y}ds, \quad dy=-\frac{\partial H_1}{\partial x}dt-\frac{\partial H_2}{\partial x}ds,\\
&dz=\frac{\partial H_1}{\partial w}dt+\frac{\partial H_2}{\partial w}ds, \quad dw=-\frac{\partial H_1}{\partial z}dt-\frac{\partial H_2}{\partial z}ds
\end{split}
\end{align}
with the symmetric Hamiltonians $H_i \in {\Bbb C}(t,s)[x,y,z,w] \ (i=1,2)$
\begin{align}\label{SdeG3H}
\begin{split}
H_1 &=H_{IV}(x,y,t;\alpha_1,\alpha_2)+\frac{\alpha_3}{t-s}xy+\frac{\alpha_2}{t-s}zw\\
&-\frac{x^2yw-2(t-s)yzw-2xyzw+yz^2w+\alpha_3yz+\alpha_2xw}{t-s},\\
H_2 &=\pi(H_1) \quad (\alpha_1+\alpha_2+\alpha_3+\alpha_4=1)
\end{split}
\end{align}
as $\varepsilon \rightarrow 0$.
\end{theorem}
Here, for notational convenience, we have renamed $X,Y,Z,W,T,S,A_1,A_2,A_3$ to $x,y,z,\\
w,t,s,\alpha_1,\alpha_2,\alpha_3$ (which are not the same as the previous $x,y,z,w,t,s,\alpha_1,\alpha_2,\alpha_3$). The transformation $\pi$ is explicitly given by
\begin{equation}
\pi:(x,y,z,w,t,s;\alpha_1,\alpha_2,\alpha_3,\alpha_4)
 \rightarrow (z,w,x,y,s,t;\alpha_1,\alpha_3,\alpha_2,\alpha_4),
\end{equation}
and the symbol $H_{IV}(x,y,t;\alpha_1,\alpha_2)$ denotes the fourth Painlev\'e Hamiltonian given by
\begin{align}
\begin{split}
&H_{IV}(x,y,t;\alpha_1,\alpha_2)=-x^2y+2xy^2-2txy-2\alpha_1y-\alpha_2x \quad (\alpha_0+\alpha_1+\alpha_2=1).
\end{split}
\end{align}

\begin{theorem}
Let us consider a polynomial Hamiltonian system with Hamiltonian\\
$H_i \in {\Bbb C}(t,s)[x,y,z,w] \ (i=1,2)$. We assume that

$(A1)$ $deg(H_i)=5$ with respect to $x,y,z,w$.

$(A2)$ This system becomes again a polynomial Hamiltonian system in each coordinate system $\{U_j,(x_j,y_j,z_j,w_j)\} \ (j=1,2,3,4)${\rm : \rm}
$$
{U_j}={{\Bbb C}^4} \ni {(x_j,y_j,z_j,w_j)} \ (j=1,2,3,4),
$$
via the following birational and symplectic transformations
\begin{align*}
&1)x_1=-(xy+zw-\alpha_1)y, \quad y_1=\frac{1}{y}, \quad z_1=zy, \quad w_1=\frac{w}{y},\\
&2)x_2=\frac{1}{x},\quad y_2=-(yx+\alpha_2)x,\quad z_2=z,\quad w_2=w,\\
&3)x_3=x, \quad y_3=y, \quad z_3=\frac{1}{z}, \quad w_3=-(zw+\alpha_3)z,\\
&4)x_4=-((x-2y-2w+2t)y+(z-2y-2w+2s)w-\alpha_4)y, \quad y_4=\frac{1}{y},\\
&z_4=(z-2y-2w+2s)y, \quad w_4=\frac{w}{y}.
\end{align*}
Then such a system coincides with the system \eqref{SdeG3}.
\end{theorem}

\begin{theorem}
The system \eqref{SdeG3} is invariant under the following transformations\rm{:\rm} with the notation $(*)=(x,y,z,w,t,s;\alpha_1,\alpha_2,\alpha_3,\alpha_4),$
\begin{align*}
        s_2: (*) &\rightarrow \left(x+\frac{\alpha_2}{y},y,z,w,t,s;\alpha_1+\alpha_2,-\alpha_2,\alpha_3,\alpha_4+\alpha_2 \right),\\
        s_3: (*) &\rightarrow \left(x,y,z+\frac{\alpha_3}{w},w,t,s;\alpha_1+\alpha_3,\alpha_2,-\alpha_3,\alpha_4+\alpha_3 \right),\\
        \pi_1: (*) &\rightarrow (z,w,x,y,s,t;\alpha_1,\alpha_3,\alpha_2,\alpha_4),\\
        \pi_2: (*) &\rightarrow (\sqrt{-1}(x-2y-2w+2t),-\sqrt{-1}y,\sqrt{-1}(z-2y-2w+2s),-\sqrt{-1}w,\\
        &-\sqrt{-1}t,-\sqrt{-1}s;\alpha_4,\alpha_2,\alpha_3,\alpha_1),\\
        \pi_3: (*) &\rightarrow (\frac{2\sqrt{-1}(xy+zw-\alpha_1)}{x},\frac{\sqrt{-1}x(xy+\alpha_2)}{2(xy+zw-\alpha_1)},\frac{2\sqrt{-1}(xy+zw-\alpha_1)}{z},\\
        &\frac{\sqrt{-1}z(zw+\alpha_3)}{2(xy+zw-\alpha_1)},-\sqrt{-1}t,-\sqrt{-1}s;-\alpha_1-\alpha_2-\alpha_3,\alpha_2,\alpha_3,1+\alpha_1).
\end{align*}
\end{theorem}

\section{A generalization of the system \eqref{SdeG3} to three variables}
In this section, we present a generalization of the system \eqref{SdeG3} to three variables, which is equivalent to the polynomial Hamiltonian system
\begin{align}\label{SdeGH3}
\begin{split}
&dx=\frac{\partial H_1}{\partial y}dt+\frac{\partial H_2}{\partial y}ds+\frac{\partial H_3}{\partial y}du, \quad dy=-\frac{\partial H_1}{\partial x}dt-\frac{\partial H_2}{\partial x}ds-\frac{\partial H_3}{\partial x}du,\\
&dz=\frac{\partial H_1}{\partial w}dt+\frac{\partial H_2}{\partial w}ds+\frac{\partial H_3}{\partial w}du, \quad dw=-\frac{\partial H_1}{\partial z}dt-\frac{\partial H_2}{\partial z}ds-\frac{\partial H_3}{\partial z}du,\\
&dq=\frac{\partial H_1}{\partial p}dt+\frac{\partial H_2}{\partial p}ds+\frac{\partial H_3}{\partial p}du, \quad dp=-\frac{\partial H_1}{\partial q}dt-\frac{\partial H_2}{\partial q}ds-\frac{\partial H_3}{\partial q}du\\
\end{split}
\end{align}
with the symmetric Hamiltonians $H_i \in {\Bbb C}(t,s,u)[x,y,z,w,q,p] \ (i=1,2,3)$
\begin{align}\label{SdeGH3H}
\begin{split}
H_1 &=H_{IV}(x,y,t;\alpha_1,\alpha_2)\\
&+R(x,y,z,w,t,s;\alpha_2,\alpha_3)+R(x,y,q,p,t,u;\alpha_2,\alpha_4),\\
H_2 &=\pi(H_1), \quad H_3=(\pi \circ \pi)(H_1) \quad (\alpha_1+\alpha_2+ \dots +\alpha_5=1),
\end{split}
\end{align}
where the transformation $\pi$ is explicitly given by
\begin{equation}
\pi:(*) \rightarrow (z,w,q,p,x,y,s,u,t;\alpha_1,\alpha_3,\alpha_4,\alpha_2,\alpha_5).
\end{equation}
Here, the symbol $(*)$ denotes $(*):=(x,y,z,w,q,p,t,s,u;\alpha_1,\alpha_2,\alpha_3,\alpha_4,\alpha_5)$, and the symbol\\
$R(q_l,p_l,q_m,p_m,t_l,t_m;\alpha,\beta)$ is explicitly given by
\begin{align*}
&R(q_l,p_l,q_m,p_m,t_l,t_m;\alpha,\beta)\\
&=\frac{\beta q_lp_l}{t_l-t_m}+\frac{\alpha q_mp_m}{t_l-t_m}\\
&-\frac{{q_l}^2p_lp_m+2t_m p_lq_mp_m-2t_l p_lq_mp_m-2q_lp_lq_mp_m+p_lq_m^2p_m+\beta p_lq_m+\alpha q_lp_m}{t_l-t_m}.
\end{align*}

\begin{theorem}
Let us consider a polynomial Hamiltonian system with Hamiltonian\\
$H_i \in {\Bbb C}(t,s,u)[x,y,z,w,q,p] \ (i=1,2,3)$. We assume that

$(A1)$ $deg(H_i)=5$ with respect to $x,y,z,w,q,p$.

$(A2)$ This system becomes again a polynomial Hamiltonian system in each coordinate system $\{U_j,(x_j,y_j,z_j,w_j,q_j,p_j)\} \ (j=1,2,..,5)${\rm : \rm}
$$
{U_j}={{\Bbb C}^6} \ni {(x_j,y_j,z_j,w_j,q_j,p_j)} \ (j=1,2,..,5),
$$
via the following birational and symplectic transformations
\begin{align*}
&1)x_1=-(yx+wz+pq-\alpha_1)y,\quad y_1=\frac{1}{y},\quad z_1=zy,\\
&w_1=\frac{w}{y}, \quad q_1=qy, \quad p_1=\frac{p}{y},\\
&2)x_2=\frac{1}{x}, \quad y_2=-(yx+\alpha_2)x, \quad z_2=z, \quad w_2=w, \quad q_2=q, \quad p_2=p,\\
&3)x_3=x, \quad y_3=y, \quad z_3=\frac{1}{z}, \quad w_3=-(wz+\alpha_3)z, \quad q_3=q, \quad p_3=p,\\
&4)x_4=x, \quad y_4=y, \quad z_4=z, \quad w_4=w, \quad q_w=\frac{1}{q}, \quad p_4=-(qp+\alpha_4)q,\\
&5)x_5=-\{(x-2y-2w-2p+2t)y+(z-2y-2w-2p+2s)w\\
&+(q-2y-2w-2p+2u)p-\alpha_5\}y, \quad y_5=\frac{1}{y}, \quad z_5=(z-2y-2w-2p+2s)y,\\
&w_5=\frac{w}{y}, \quad q_5=(q-2y-2w-2p+2u)y, \quad p_5=\frac{p}{y}.
\end{align*}
Then such a system coincides with the system \eqref{SdeGH3}.
\end{theorem}

\begin{theorem}
The system \eqref{SdeGH3} is invariant under the following transformations\rm{:\rm} with the notation $(*)=(x,y,z,w,q,p,t,s,u;\alpha_1,\alpha_2,\ldots,\alpha_5),$
\begin{align*}
        s_2: (*) &\rightarrow \left(x+\frac{\alpha_2}{y},y,z,w,q,p,t,s,u;\alpha_1+\alpha_1,-\alpha_2,\alpha_3,\alpha_4,\alpha_5+\alpha_2 \right),\\
        s_3: (*) &\rightarrow \left(x,y,z+\frac{\alpha_3}{w},w,t,s,u;\alpha_1+\alpha_3,\alpha_2,-\alpha_3,\alpha_4,\alpha_5+\alpha_3 \right),\\
        s_4: (*) &\rightarrow \left(x,y,z,w,q+\frac{\alpha_4}{p},p,t,s,u;\alpha_1+\alpha_4,\alpha_2,\alpha_3,-\alpha_4,\alpha_5+\alpha_4 \right),\\
        \pi_1: (*) &\rightarrow (z,w,q,p,x,y,s,u,t;\alpha_1,\alpha_3,\alpha_4,\alpha_2,\alpha_5),\\
        \pi_2: (*) &\rightarrow (\sqrt{-1}(x-2y-2w-2p+2t),-\sqrt{-1}y,\sqrt{-1}(z-2y-2w-2p+2s),\\
        &-\sqrt{-1}w,\sqrt{-1}(q-2y-2w-2p+2u),-\sqrt{-1}p,-\sqrt{-1}t,-\sqrt{-1}s,\\
        &-\sqrt{-1}u;\alpha_5,\alpha_2,\alpha_3,\alpha_4,\alpha_1),\\
        \pi_3: (*) &\rightarrow (\frac{2\sqrt{-1}(xy+zw+qp-\alpha_1)}{x},\frac{\sqrt{-1}x(xy+\alpha_2)}{2(xy+zw+qp-\alpha_1)},\\
        &\frac{2\sqrt{-1}(xy+zw+qp-\alpha_1)}{z},\frac{\sqrt{-1}z(zw+\alpha_3)}{2(xy+zw+qp-\alpha_1)},\\
        &\frac{2\sqrt{-1}(xy+zw+qp-\alpha_1)}{q},\frac{\sqrt{-1}q(qp+\alpha_4)}{2(xy+zw+qp-\alpha_1)},\\
        &-\sqrt{-1}t,-\sqrt{-1}s,-\sqrt{-1}u;-\alpha_1-\alpha_2-\alpha_3-\alpha_4,\alpha_2,\alpha_3,\alpha_4,1+\alpha_1).
        \end{align*}
\end{theorem}

\section{Another degeneration from the system \eqref{dVV}}

As the fourth-order analogue of the above confluence process from $P_{V}$ to $P_{III}$ (see \cite{T2}), we consider the following coupling confluence process from the system \eqref{dVV}. We take the following coupling confluence process $P_{V} \rightarrow P_{III}$ for each coordinate system $(x,y)$ and $(z,w)$ of the system \eqref{dVV}.

\begin{theorem}
For the system \eqref{dVV}, we make the change of parameters and variables
\begin{gather}
\begin{gathered}
\alpha_1=A_0, \quad \alpha_2=A_2, \quad \alpha_3=\frac{1}{\varepsilon}, \quad \alpha_4=2A_1-\frac{1}{\varepsilon}, \quad \alpha_5=A_3
\end{gathered}\\
\begin{gathered}
t=-\varepsilon T, \quad s=-\varepsilon S, \quad X=-t(x-1), \quad Z=-s(z-1),\\
Y=-\frac{y}{t}, \quad W=-\frac{w}{s}
\end{gathered}
\end{gather}
from $\alpha_1,\alpha_2, \dots,\alpha_5,t,s,x,y,z,w$ to $A_0,\dots ,A_3,\varepsilon,T,S,X,Y,Z,W$. Then this system can also be written in the new variables $T,S,X,Y,Z,W$ and parameters $A_0,\dots,A_3,\varepsilon$ as a Hamiltonian system. This new system tends to the polynomial Hamiltonian system
\begin{align}\label{SdeGa}
\begin{split}
&dx=\frac{\partial H_1}{\partial y}dt+\frac{\partial H_2}{\partial y}ds, \quad dy=-\frac{\partial H_1}{\partial x}dt-\frac{\partial H_2}{\partial x}ds,\\
&dz=\frac{\partial H_1}{\partial w}dt+\frac{\partial H_2}{\partial w}ds, \quad dw=-\frac{\partial H_1}{\partial z}dt-\frac{\partial H_2}{\partial z}ds
\end{split}
\end{align}
with the symmetric Hamiltonians $H_i \in {\Bbb C}(t,s)[x,y,z,w] \ (i=1,2)$
\begin{align}\label{SdeGaH4545}
\begin{split}
H_1 &=H_{III}(x,y,t;\alpha_0,\alpha_1)+\frac{\alpha_2 s}{t(t-s)}xy+\frac{\alpha_0}{t-s}zw\\
&-\frac{tyz^2w+sx^2yw-2txyzw+\alpha_2tyz+\alpha_0sxw}{t(t-s)},\\
H_2 &=\pi(H_1) \quad (\alpha_0+2\alpha_1+\alpha_2+\alpha_3=1),
\end{split}
\end{align}
as $\varepsilon \rightarrow 0$.
\end{theorem}
Here, for notational convenience, we have renamed $X,Y,Z,W,T,S,A_1,A_2,A_3$ to $x,y,z,\\
w,t,s,\alpha_1,\alpha_2,\alpha_3$ (which are not the same as the previous $x,y,z,w,t,s,\alpha_1,\alpha_2,\alpha_3$). The transformation $\pi$ is explicitly given by
\begin{equation}
\pi:(x,y,z,w,t,s;\alpha_0,\alpha_1,\alpha_2,\alpha_3)
 \rightarrow (z,w,x,y,s,t;\alpha_2,\alpha_1,\alpha_0,\alpha_3),
\end{equation}
and the symbol $H_{III}(x,y,t;\alpha_0,\alpha_1,\alpha_2)$ denotes the third Painlev\'e Hamiltonian given by
\begin{align}
\begin{split}
&H_{III}(x,y,t;\alpha_0,\alpha_1)=\frac{x^2y(y-1)+x\{(1-2\alpha_1)y-\alpha_0\}+ty}{t} \quad (\alpha_0+2\alpha_1+\alpha_2=1).
\end{split}
\end{align}

\begin{theorem}
Let us consider a polynomial Hamiltonian system with Hamiltonian\\
$H_i \in {\Bbb C}(t,s)[x,y,z,w] \ (i=1,2)$. We assume that

$(A1)$ $deg(H_i)=5$ with respect to $x,y,z,w$.

$(A2)$ This system becomes again a polynomial Hamiltonian system in each coordinate system $\{U_j,(x_j,y_j,z_j,w_j)\} \ (j=0,1,2,3)${\rm : \rm}
$$
{U_j}={{\Bbb C}^4} \ni {(x_j,y_j,z_j,w_j)} \ (j=0,1,2,3),
$$
via the following birational and symplectic transformations
\begin{align*}
&0)x_0=\frac{1}{x}, \quad y_0=-(yx+\alpha_0)x, \quad z_0=z, \quad w_0=w,\\
&1)x_1=x, \quad y_1=y+\frac{s}{t}w+\frac{2\{(z-\frac{s}{t}x)w-\alpha_1\}}{x}+\frac{t}{x^2},\\
&z_1=\frac{z-\frac{s}{t}x}{x^2}, \quad w_1=x^2w,\\
&2)x_2=x,\quad y_2=y,\quad z_2=\frac{1}{z},\quad w_2=-(zw+\alpha_2)z,\\
&3)x_3=\frac{1}{x}, \quad y_3=-((y+w-1)x+\alpha_3)x, \quad z_3=z-x, \quad w_3=w.
\end{align*}
Then such a system coincides with the system \eqref{SdeGa}.
\end{theorem}

\begin{theorem}
The system \eqref{SdeGa} is invariant under the following transformations\rm{:\rm} with the notation $(*)=(x,y,z,w,t,s;\alpha_1,\alpha_2,\alpha_3,\alpha_4),$
\begin{align*}
        s_0: (*) &\rightarrow \left(x+\frac{\alpha_0}{y},y,z,w,t,s;-\alpha_0,\alpha_1+\alpha_0,\alpha_2,\alpha_3 \right),\\
        s_2: (*) &\rightarrow \left(x,y,z+\frac{\alpha_2}{w},w,t,s;\alpha_0,\alpha_1+\alpha_2,-\alpha_2,\alpha_3 \right),\\
        s_3: (*) &\rightarrow \left(x+\frac{\alpha_3}{y+w-1},y,z+\frac{\alpha_3}{y+w-1},w,t,s;\alpha_0,\alpha_1+\alpha_3,\alpha_2,-\alpha_3 \right),\\
        \pi_1: (*) &\rightarrow (z,w,x,y,s,t;\alpha_2,\alpha_1,\alpha_0,\alpha_3),\\
        \pi_2: (*) &\rightarrow \left(\frac{t}{x},-\frac{(xy+\alpha_0)x}{t},\frac{s}{z},-\frac{(zw+\alpha_2)z}{s},t,s;\alpha_0,\alpha_1+\alpha_3-\frac{1}{2},\alpha_2,1-\alpha_3 \right),\\
        \pi_3: (*) &\rightarrow (x-z,y,-z,1-y-w,t-s,-s;\alpha_0,\alpha_1,\alpha_3,\alpha_2).
\end{align*}
\end{theorem}

\section{A generalization of the system \eqref{SdeGa} to three variables}
In this section, we present a generalization of the system \eqref{SdeGa} to three variables, which is equivalent to the polynomial Hamiltonian system
\begin{align}
\begin{split}
&dx=\frac{\partial H_1}{\partial y}dt+\frac{\partial H_2}{\partial y}ds+\frac{\partial H_3}{\partial y}du, \quad dy=-\frac{\partial H_1}{\partial x}dt-\frac{\partial H_2}{\partial x}ds-\frac{\partial H_3}{\partial x}du,\\
&dz=\frac{\partial H_1}{\partial w}dt+\frac{\partial H_2}{\partial w}ds+\frac{\partial H_3}{\partial w}du, \quad dw=-\frac{\partial H_1}{\partial z}dt-\frac{\partial H_2}{\partial z}ds-\frac{\partial H_3}{\partial z}du,\\
&dq=\frac{\partial H_1}{\partial p}dt+\frac{\partial H_2}{\partial p}ds+\frac{\partial H_3}{\partial p}du, \quad dp=-\frac{\partial H_1}{\partial q}dt-\frac{\partial H_2}{\partial q}ds-\frac{\partial H_3}{\partial q}du\\
\end{split}
\end{align}
with the symmetric Hamiltonians $H_i \in {\Bbb C}(t,s,u)[x,y,z,w,q,p] \ (i=1,2,3)$
\begin{align}\label{SdeGaH}
\begin{split}
H_1 &=H_{III}(x,y,t;\alpha_0,\alpha_1)\\
&+R(x,y,z,w,t,s;\alpha_0,\alpha_2)+R(x,y,q,p,t,u;\alpha_0,\alpha_4),\\
H_2 &=\pi(H_1), \quad H_3=(\pi \circ \pi)(H_1) \quad (\alpha_0+2\alpha_1+\alpha_2+\alpha_3+\alpha_4=1),
\end{split}
\end{align}
where the transformation $\pi$ is explicitly given by
\begin{equation}
\pi:(*) \rightarrow (z,w,q,p,x,y,s,u,t;\alpha_2,\alpha_1,\alpha_4,\alpha_3,\alpha_0).
\end{equation}
Here the symbol $(*)$ denotes $(*):=(x,y,z,w,q,p,t,s,u;\alpha_0,\alpha_1,\alpha_2,\alpha_3,\alpha_4)$, and the symbol\\
$R(q_l,p_l,q_m,p_m,t_l,t_m;\alpha,\beta)$ is explicitly given by
\begin{align*}
&R(q_l,p_l,q_m,p_m,t_l,t_m;\alpha,\beta)\\
&=\frac{\beta t_m q_lp_l}{t_l(t_l-t_m)}+\frac{\alpha q_mp_m}{t_l-t_m}\\
&-\frac{t_l p_lq_m^2p_m+t_mq_l^2p_lp_m-2t_lq_lp_lq_mp_m+\beta t_lp_lq_m+\alpha t_mq_lp_m}{t_l(t_l-t_m)}.
\end{align*}

\begin{theorem}
Let us consider a polynomial Hamiltonian system with Hamiltonian\\
$H_i \in {\Bbb C}(t,s,u)[x,y,z,w,q,p] \ (i=1,2,3)$. We assume that

$(A1)$ $deg(H_i)=5$ with respect to $x,y,z,w,q,p$.

$(A2)$ This system becomes again a polynomial Hamiltonian system in each coordinate system $\{U_j,(x_j,y_j,z_j,w_j,q_j,p_j)\} \ (j=0,1,..,4)${\rm : \rm}
$$
{U_j}={{\Bbb C}^6} \ni {(x_j,y_j,z_j,w_j,q_j,p_j)} \ (j=0,1,..,4),
$$
via the following birational and symplectic transformations
\begin{align*}
&0)x_0=\frac{1}{x}, \quad y_0=-(yx+\alpha_0)x, \quad z_0=z, \quad w_0=w, \quad q_0=q, \quad p_0=p,\\
&1)x_1=x, \quad y_1=y+\frac{s}{t}w+\frac{u}{t}p+\frac{2\{(z-\frac{s}{t}x)w+(q-\frac{u}{t}x)p-\alpha_1\}}{x}+\frac{t}{x^2},\\
&z_1=\frac{z-\frac{s}{t}x}{x^2}, \quad w_1=x^2w, \quad q_1=\frac{q-\frac{u}{t}x}{x^2}, \quad p_1=x^2p,\\
&2)x_2=x,\quad y_2=y,\quad z_2=\frac{1}{z},\quad w_2=-(zw+\alpha_2)z, \quad q_2=q, \quad p_2=p,\\
&3)x_3=\frac{1}{x}, \quad y_3=-((y+w+p-1)x+\alpha_3)x, \quad z_3=z-x,\\
&w_3=w. \quad q_3=q-x, \quad p_3=p,\\
&4)x_4=x, \quad y_4=y, \quad z_4=z, \quad w_4=w, \quad q_4=\frac{1}{q}, \quad p_4=-(qp+\alpha_4)q.
\end{align*}
Then such a system coincides with the system \eqref{SdeGaH}.
\end{theorem}

\begin{theorem}
The system \eqref{SdeGaH} is invariant under the following transformations\rm{:\rm} with the notation $(*)=(x,y,z,w,q,p,t,s,u;\alpha_0,\alpha_1,\ldots,\alpha_4),$
\begin{align*}
        s_0: (*) &\rightarrow \left(x+\frac{\alpha_0}{y},y,z,w,q,p,t,s,u;-\alpha_0,\alpha_1+\alpha_0,\alpha_2,\alpha_3,\alpha_4 \right),\\
        s_2: (*) &\rightarrow \left(x,y,z+\frac{\alpha_2}{w},w,q,p,t,s,u;\alpha_0,\alpha_1+\alpha_2,-\alpha_2,\alpha_3,\alpha_4 \right),\\
        s_3: (*) &\rightarrow (x+\frac{\alpha_3}{y+w+p-1},y,z+\frac{\alpha_3}{y+w+p-1},w,q+\frac{\alpha_3}{y+w+p-1},p,\\
        &t,s,u;\alpha_0,\alpha_1+\alpha_3,\alpha_2,-\alpha_3,\alpha_4),\\
        s_4: (*) &\rightarrow \left(x,y,z,w,q+\frac{\alpha_4}{p},p,t,s,u;\alpha_0,\alpha_1+\alpha_4,\alpha_2,\alpha_3,-\alpha_4 \right),\\
        \pi_1: (*) &\rightarrow (z,w,q,p,x,y,s,u,t;\alpha_2,\alpha_1,\alpha_4,\alpha_3,\alpha_0),\\
        \pi_2: (*) &\rightarrow (\frac{t}{x},-\frac{(xy+\alpha_0)x}{t},\frac{s}{z},-\frac{(zw+\alpha_2)z}{s},\frac{u}{q},-\frac{(qp+\alpha_4)q}{u},\\
        &t,s,u;\alpha_0,\alpha_1+\alpha_3-\frac{1}{2},\alpha_2,1-\alpha_3,\alpha_4).
        \end{align*}
\end{theorem}

\section{Other generalization of the third Painlev\'e system}
In this section, we find a generalization of the third Painlev\'e system to two variables $t,s$, which is different from the system \eqref{SdeGa}. This system is equivalent to the Hamiltonian system
\begin{align}\label{ASdeGa}
\begin{split}
&dx=\frac{\partial H_1}{\partial y}dt+\frac{\partial H_2}{\partial y}ds, \quad dy=-\frac{\partial H_1}{\partial x}dt-\frac{\partial H_2}{\partial x}ds,\\
&dz=\frac{\partial H_1}{\partial w}dt+\frac{\partial H_2}{\partial w}ds, \quad dw=-\frac{\partial H_1}{\partial z}dt-\frac{\partial H_2}{\partial z}ds
\end{split}
\end{align}
with the polynomial Hamiltonians $H_i \in {\Bbb C}(t,s)[x,y,z,w] \ (i=1,2)$
\begin{align}\label{ASdeGaH}
\begin{split}
H_1 &=\frac{-x^3y^2+sx^2y^2-(2\alpha_1+\alpha_2)x^2y+\{(2\alpha_1+\alpha_2)s+\eta_1 t\}xy}{ts}-\frac{\alpha_1(\alpha_1+\alpha_2)x+\eta_1 tsy}{ts}\\
&+\frac{z^3w^2-tz^2w^2+(2\alpha_1+\alpha_2)z^2w+(\alpha_3 t-\eta_0 s)zw}{t^2}+\frac{\alpha_1(\alpha_1+\alpha_2)z+\eta_0 tsw}{t^2}\\
&-\frac{\eta_0 s-\alpha_3 t}{t^2}xy+\frac{\eta_1}{s}zw-\frac{xz\{2(tx-sz)yw-sxy^2+tzw^2-(2\alpha_1+\alpha_2)(sy-tw)\}}{t^2s},\\
H_2 &=\pi(H_1) \quad (2\alpha_1+\alpha_2+\alpha_3+\alpha_4=1),
\end{split}
\end{align}
where the transformation $\pi$ is explicitly given by
\begin{equation}
\pi:(x,y,z,w,t,s;\eta_0,\eta_1,\alpha_1,\alpha_2,\alpha_3,\alpha_4)
 \rightarrow (z,w,x,y,s,t;\eta_1,\eta_0,\alpha_1,\alpha_2,\alpha_4,\alpha_3).
\end{equation}

\begin{theorem}
Let us consider a polynomial Hamiltonian system with Hamiltonian\\
$H_i \in {\Bbb C}(t,s)[x,y,z,w] \ (i=1,2)$. We assume that

$(A1)$ $deg(H_i)=5$ with respect to $x,y,z,w$.

$(A2)$ This system becomes again a polynomial Hamiltonian system in each coordinate system $\{U_j,(x_j,y_j,z_j,w_j)\} \ (j=1,2,3,4)${\rm : \rm}
$$
{U_j}={{\Bbb C}^4} \ni {(x_j,y_j,z_j,w_j)} \ (j=1,2,3,4),
$$
via the following birational and symplectic transformations
\begin{align*}
&1)x_1=\frac{1}{x}, \quad y_1=-(xy+zw+\alpha_1)x, \quad z_1=\frac{z}{x}, \quad w_1=xw,\\
&2)x_2=\frac{1}{x},\quad y_2=-(xy+zw+\alpha_1+\alpha_2)x,\quad z_2=\frac{z}{x},\quad w_2=xw,\\
&3)x_3=x, \quad y_3=y-\frac{\eta_0}{z}, \quad z_3=z, \quad w_3=w-\frac{\alpha_3}{z}+\frac{\eta_0(x-s)}{z^2},\\
&4)x_4=x, \quad y_4=y-\frac{\alpha_4}{x}+\frac{\eta_1(z-t)}{x^2}, \quad z_4=z, \quad w_4=w-\frac{\eta_1}{x}.
\end{align*}
Then such a system coincides with the system \eqref{ASdeGa}.
\end{theorem}

\begin{theorem}
The system \eqref{ASdeGa} is invariant under the following transformations\rm{:\rm} with the notation $(*)=(x,y,z,w,t,s;\eta_0,\eta_1,\alpha_1,\alpha_2,\alpha_3,\alpha_4),$
\begin{align*}
        s_2: (*) &\rightarrow (x,y,z,w,t,s;\eta_0,\eta_1,\alpha_1+\alpha_2,-\alpha_2,\alpha_3,\alpha_4),\\
        s_3: (*) &\rightarrow \left(x,y-\frac{\eta_0}{z},z,w-\frac{\alpha_3}{z}+\frac{\eta_0(x-s)}{z^2},t,s;-\eta_0,\eta_1,\alpha_1+\alpha_3,\alpha_2,-\alpha_3,\alpha_4 \right),\\
        s_4: (*) &\rightarrow \left(x,y-\frac{\alpha_4}{x}+\frac{\eta_1(z-t)}{x^2},z,w-\frac{\eta_1}{x},t,s;\eta_0,-\eta_1,\alpha_1+\alpha_4,\alpha_2,\alpha_3,-\alpha_4 \right).
\end{align*}
\end{theorem}

\section{Degeneration from the system \eqref{SdeG3}}

As the fourth-order analogue of the above confluence process from $P_{IV}$ to $P_{II}$ (see \cite{T2}), we consider the following coupling confluence process from the system \eqref{SdeG3}. We take the following coupling confluence process $P_{IV} \rightarrow P_{II}$ for each coordinate system $(x,y)$ and $(z,w)$ of the system \eqref{SdeG3}.

\begin{theorem}
For the system \eqref{SdeG3}, we make the change of parameters and variables
\begin{gather}
\begin{gathered}
\alpha_1=\frac{1}{4{\varepsilon}^6}, \quad \alpha_2=A_1, \quad \alpha_3=A_3, \quad \alpha_4=A_2-\frac{1}{4{\varepsilon}^6},
\end{gathered}\\
\begin{gathered}
t=-\frac{1-{\varepsilon}^4 T}{\sqrt{2}{\varepsilon}^3}, \quad s=-\frac{1-{\varepsilon}^4 S}{\sqrt{2}{\varepsilon}^3}, \quad x=\frac{1+2{\varepsilon}^2 X}{\sqrt{2} {\varepsilon}^3},\\
z=\frac{1+2{\varepsilon}^2 Z}{\sqrt{2} {\varepsilon}^3}, \quad y=\frac{\varepsilon Y}{\sqrt{2}}, \quad w=\frac{\varepsilon W}{\sqrt{2}}
\end{gathered}
\end{gather}
from $\alpha_1,\alpha_2, \dots,\alpha_4,t,s,x,y,z,w$ to $A_1,A_2,A_3,\varepsilon,T,S,X,Y,Z,W$. Then this system can also be written in the new variables $T,S,X,Y,Z,W$ and parameters $A_1,A_2,A_3,\varepsilon$ as a Hamiltonian system. This new system tends to the polynomial Hamiltonian system
\begin{align}\label{SdeG4}
\begin{split}
&dx=\frac{\partial H_1}{\partial y}dt+\frac{\partial H_2}{\partial y}ds, \quad dy=-\frac{\partial H_1}{\partial x}dt-\frac{\partial H_2}{\partial x}ds,\\
&dz=\frac{\partial H_1}{\partial w}dt+\frac{\partial H_2}{\partial w}ds, \quad dw=-\frac{\partial H_1}{\partial z}dt-\frac{\partial H_2}{\partial z}ds
\end{split}
\end{align}
with the symmetric Hamiltonians $H_i \in {\Bbb C}(t,s)[x,y,z,w] \ (i=1,2)$
\begin{align}\label{SdeG4H}
\begin{split}
H_1 &=H_{II}(x,y,t;\alpha_3)+\frac{\alpha_1}{t-s}xy-\frac{\alpha_1}{t-s}yz-\frac{\alpha_3}{t-s}xw+\frac{\alpha_3}{t-s}zw-\frac{\{2(x-z)^2-(t-s)\}yw}{2(t-s)},\\
H_2 &=\pi(H_1) \quad (\alpha_1+\alpha_2+\alpha_3=1)
\end{split}
\end{align}
as $\varepsilon \rightarrow 0$.
\end{theorem}
Here, for notational convenience, we have renamed $X,Y,Z,W,T,S,A_1,A_2$ to $x,y,z,w,\\
t,s,\alpha_1,\alpha_2$ (which are not the same as the previous $x,y,z,w,t,s,\alpha_1,\alpha_2$). The transformation $\pi$ is explicitly given by
\begin{equation}
\pi:(x,y,z,w,t,s;\alpha_1,\alpha_2,\alpha_3)
 \rightarrow (z,w,x,y,s,t;\alpha_3,\alpha_2,\alpha_1),
\end{equation}
and the symbol $H_{II}(x,y,t;\alpha_0,\alpha_1)$ denotes the second Painlev\'e Hamiltonian given by
\begin{equation}
H_{II}(x,y,t;\alpha_1)=\frac{1}{2}y^2-\left(x^2+\frac{t}{2}\right)y-\alpha_1 x.
\end{equation}

\begin{theorem}
Let us consider a polynomial Hamiltonian system with Hamiltonian\\
$H_i \in {\Bbb C}(t,s)[x,y,z,w] \ (i=1,2)$. We assume that

$(A1)$ $deg(H_i)=5$ with respect to $x,y,z,w$.

$(A2)$ This system becomes again a polynomial Hamiltonian system in each coordinate system $\{U_j,(x_j,y_j,z_j,w_j)\} \ (j=1,2,3)${\rm : \rm}
$$
{U_j}={{\Bbb C}^4} \ni {(x_j,y_j,z_j,w_j)} \ (j=1,2,3),
$$
via the following birational and symplectic transformations
\begin{align*}
&1)x_1=\frac{1}{x}, \quad y_1=-(yx+\alpha_3)x, \quad z_1=z, \quad w_1=w,\\
&2)x_2=\frac{1}{x}, \quad y_2=-\{(y+w-2x^2-t)x-2\left((z-x)x-\frac{t-s}{4}\right)\left(\frac{w}{x}\right)+\alpha_2\}x,\\
&z_2=\left((z-x)x-\frac{t-s}{2}\right)x, \quad w_2=\frac{w}{x^2},\\
&3)x_3=x,\quad y_3=y,\quad z_3=\frac{1}{z},\quad w_3=-(zw+\alpha_1)z.
\end{align*}
Then such a system coincides with the system \eqref{SdeG4}.
\end{theorem}

\begin{theorem}
The system \eqref{SdeG4} is invariant under the following transformations\rm{:\rm} with the notation $(*)=(x,y,z,w,t,s;\alpha_1,\alpha_2,\alpha_3,\alpha_4),$
\begin{align*}
        s_1: (*) &\rightarrow \left(x+\frac{\alpha_3}{y},y,z,w,t,s;\alpha_1,\alpha_2+\alpha_3,-\alpha_3 \right),\\
        s_3: (*) &\rightarrow \left(x,y,z+\frac{\alpha_1}{w},w,t,s;-\alpha_1,\alpha_2+\alpha_1,\alpha_3 \right),\\
        \pi_1: (*) &\rightarrow (z,w,x,y,s,t;\alpha_3,\alpha_2,\alpha_1).
\end{align*}
\end{theorem}

For the system \eqref{SdeG4}, we make the change of variables
\begin{equation}
t=T, \quad s=T+S
\end{equation}
from $t,s,x,y,z,w$ to $T,S,x,y,z,w$. Then this system can also be written in the new variables $T,S,x,y,z,w$ as the Hamiltonian system
\begin{align}
\begin{split}
dx=&\left(-x^2+y+w-\frac{T}{2} \right)dT+\left(-\frac{(x-z)(xw-zw-\alpha_1)}{S}+\frac{w}{2} \right)dS,\\
dy=&(2xy+\alpha_3)dT+\left(\frac{2xyw-2yzw-\alpha_1 y+\alpha_3 w}{S} \right)dS,\\
dz=&\left(-z^2+w+y-\frac{S}{2}-\frac{T}{2} \right)dT\\
&+\left(-z^2+w-\frac{S}{2}-\frac{T}{2}+\frac{y}{2}-\frac{(X-Z)(XY-YZ+\alpha_3)}{S} \right)dS,\\
dw=&(2zw+\alpha_1)dT+\left(2zw+\alpha_1-\frac{2xyw-2yzw-\alpha_1 y+\alpha_3 w}{S} \right)dS
\end{split}
\end{align}
with the polynomial Hamiltonians
\begin{align}\label{Kimura}
H_1&=-x^2y+\frac{y^2}{2}-\frac{Ty}{2}-\alpha_3 x-z^2w+\frac{w^2}{2}-\frac{Sw}{2}-\alpha_1 z-\frac{Tw}{2}+yw,\\
\begin{split}
H_2&=-\frac{(x-z)(xw-zw-\alpha_1)y}{S}+\frac{yw}{2}-\frac{Tw}{2}-\frac{\alpha_3 xw}{S}+\frac{\alpha_3 zw}{S}\\
&-z^2w+\frac{w^2}{2}-\frac{Sw}{2}-\alpha_1 z.
\end{split}
\end{align}

\section{Autonomous version of the system \eqref{SdeG4}}

In this section, we present an autonomous version of the system  \eqref{SdeG4} given by
\begin{equation}\label{eq:autoG(1,4)}
  \left\{
  \begin{aligned}
   dq_1 =&\frac{\partial K_1}{\partial p_1}dt+\frac{\partial K_2}{\partial p_1}ds,\\
   dp_1 =&-\frac{\partial K_1}{\partial q_1}dt-\frac{\partial K_2}{\partial q_1}ds,\\
   dq_2 =&\frac{\partial K_1}{\partial p_2}dt+\frac{\partial K_2}{\partial p_2}ds,\\
   dp_2 =&-\frac{\partial K_1}{\partial q_2}dt-\frac{\partial K_2}{\partial q_2}ds
   \end{aligned}
  \right. 
\end{equation}
with the polynomial Hamiltonians
\begin{align}
\begin{split}
K_1=&-q_1^2 p_1+\frac{p_1^2}{2}-\alpha_2 q_1-q_2^2 p_2+\frac{p_2^2}{2}-\alpha_0 q_2+p_1 p_2,\\
K_2=&q_1^2 p_1 p_2+q_2^2 p_1 p_2-2q_1p_1q_2p_2-\alpha_0 q_1p_1-\alpha_2 q_2p_2+\alpha_0 p_1q_2+\alpha_2 q_1p_2.
\end{split}
\end{align}

\begin{proposition}
The system \eqref{eq:autoG(1,4)} satisfies the compatibility conditions$:$
\begin{equation}
\frac{\partial }{\partial s} \frac{\partial q_1}{\partial t}=\frac{\partial }{\partial t} \frac{\partial q_1}{\partial s}, \quad \frac{\partial }{\partial s} \frac{\partial p_1}{\partial t}=\frac{\partial }{\partial t} \frac{\partial p_1}{\partial s}, \quad \frac{\partial }{\partial s} \frac{\partial q_2}{\partial t}=\frac{\partial }{\partial t} \frac{\partial q_2}{\partial s}, \quad \frac{\partial }{\partial s} \frac{\partial p_2}{\partial t}=\frac{\partial }{\partial t} \frac{\partial p_2}{\partial s}.
\end{equation}
\end{proposition}

\begin{proposition}
The system \eqref{eq:autoG(1,4)} has $K_1$ and $K_2$ as its first integrals.
\end{proposition}

\begin{proposition}
Two Hamiltonians $K_1$ and $K_2$ satisfy
\begin{equation}
\{K_1,K_2\}=0,
\end{equation}
where
\begin{equation}
\{K_1,K_2\}=\frac{\partial K_1}{\partial p_1}\frac{\partial K_2}{\partial q_1}-\frac{\partial K_1}{\partial q_1}\frac{\partial K_2}{\partial p_1}+\frac{\partial K_1}{\partial p_2}\frac{\partial K_2}{\partial q_2}-\frac{\partial K_1}{\partial q_2}\frac{\partial K_2}{\partial p_2}.
\end{equation}
\end{proposition}
Here, $\{,\}$ denotes the poisson bracket such that $\{p_i,q_j\}={\delta}_{ij}$ (${\delta}_{ij}$:kronecker's delta).

\begin{theorem}
The system \eqref{eq:autoG(1,4)} admits the extended affine Weyl group symmetry of type $D_3^{(2)}$ as the group of its B{\"a}cklund transformations, whose generators $s_0,s_1,s_2,{\pi}$ defined as follows$:$ with {\it the notation} $(*):=(q_1,p_1,q_2,p_2,t,s;\alpha_0,\alpha_1,\alpha_2)$\rm{: \rm}
\begin{align*}
s_0:(*) \rightarrow &\left(q_1,p_1,q_2+\frac{\alpha_0}{p_2},p_2,t,s;-\alpha_0,\alpha_1+2\alpha_0,\alpha_2 \right),\\
s_1:(*) \rightarrow &\left(q_1,p_1-\frac{\alpha_1}{q_1-q_2},q_2,p_2+\frac{\alpha_1}{q_1-q_2},t,s;\alpha_0+\alpha_1,-\alpha_1,\alpha_2+\alpha_1 \right),\\
s_2:(*) \rightarrow &\left(q_1+\frac{\alpha_2}{p_1},p_1,q_2,p_2,t,s;\alpha_0,\alpha_1+2\alpha_2,-\alpha_2 \right),\\
\pi:(*) \rightarrow &(q_2,p_2,q_1,p_1,t,s;\alpha_2,\alpha_1,\alpha_0),
\end{align*}
where the parameters $\alpha_i$ satisfy the relation $\alpha_0+\alpha_1+\alpha_2=0$.
\end{theorem}

\begin{theorem}
Let us consider a polynomial Hamiltonian system with Hamiltonian $K \in {\Bbb C}[q_1,p_1,q_2,p_2]$. We assume that

$(C1)$ $deg(K)=4$ with respect to $q_1,p_1,q_2,p_2$.

$(C2)$ This system becomes again a polynomial Hamiltonian system in each coordinate $R_i \ (i=0,1,2)${\rm : \rm}
\begin{align*}
\begin{split}
R_0:(x_0,y_0,z_0,w_0)=&\left(q_1,p_1,\frac{1}{q_2},-(q_2p_2+\alpha_0)q_2 \right),\\
R_1:(x_1,y_1,z_1,w_1)=&\left(\frac{1}{q_1},-((p_1+p_2-2q_1^2)q_1-2(q_2-q_1)p_2+\alpha_1)q_1,(q_2-q_1)q_1^2,\frac{p_2}{q_1^2} \right),\\
R_2:(x_2,y_2,z_2,w_2)=&\left(\frac{1}{q_1},-(q_1p_1+\alpha_2)q_1,q_2,p_2 \right),
\end{split}
\end{align*}
where the parameters $\alpha_i$ satisfy the relation $\alpha_0+\alpha_1+\alpha_2=0$. Then such a system coincides with the Hamiltonian system \eqref{eq:autoG(1,4)} with the polynomial Hamiltonians $K_1,K_2$.
\end{theorem}
We note that the conditions $(C2)$ should be read that
\begin{align*}
\begin{split}
&R_i \ (i=0,1,2)
\end{split}
\end{align*}
are polynomials with respect to $x_i,y_i,z_i,w_i$.

Next, let us consider the relation between the polynomial Hamiltonian system \eqref{eq:autoG(1,4)} and soliton equations. In this paper, we can make the birational transformations between the polynomial Hamiltonian system \eqref{eq:autoG(1,4)} and soliton equations.
\begin{theorem}
The birational transformations
\begin{equation}\label{eq:G(1,4)2}
  \left\{
  \begin{aligned}
   x =&q_1,\\
   y =&p_1+p_2-q_1^2,\\
   z =&2q_1^3-2q_1p_2+2q_2p_2+\alpha_0+\alpha_2,\\
   w =&-6q_1^4+8q_1^2 p_2+6q_1^2 p_1+2q_2^2 p_2-4q_1q_2p_2-2\alpha_0(q_1-q_2),\\
   S =&\frac{1}{2}s
   \end{aligned}
  \right. 
\end{equation}
take the Hamiltonian system \eqref{eq:autoG(1,4)} to the system
\begin{equation}\label{eq:G(1,4)3}
  \left\{
  \begin{aligned}
   dx =&y dt+(w-6x^2y)dS,\\
   dy =&z dt+(-2x(w-6x^2y))dS,\\
   dz =&w dt+(2(2x^2-y)(w-6x^2y))dS,\\
   dw =&(12x^3y+12xy^2+6x^2z-2xw)dt+(-2(4x^3-6xy+z)(w-6x^2y))dS.
   \end{aligned}
  \right. 
\end{equation}
\end{theorem}
Setting $u:=x$, we see that
\begin{equation}
\frac{\partial u}{\partial t}=y, \quad  \frac{\partial^2 u}{\partial t^2}=z, \quad \frac{\partial^3 u}{\partial t^3}=w,
\end{equation}
and
\begin{equation}\label{eq:G(1,4)4}
  \left\{
  \begin{aligned}
   \frac{\partial^4 u}{\partial t^4} =&12u^3 \frac{\partial u}{\partial t}+12u\left(\frac{\partial u}{\partial t} \right)^2+6u^2 \frac{\partial^2 u}{\partial t^2}-2u \frac{\partial^3 u}{\partial t^3},\\
   \frac{\partial u}{\partial S} =&\frac{\partial^3 u}{\partial t^3}-6u^2\frac{\partial u}{\partial t},
   \end{aligned}
  \right. 
\end{equation}
where the second equation just coincides with the mKdV equation.

Integrating for the first equation of \eqref{eq:G(1,4)4}, we obtain
\begin{equation}\label{eq:G(1,4)5}
  \left\{
  \begin{aligned}
   \frac{\partial^3 u}{\partial t^3} =&3u^4+6u^2 \frac{\partial u}{\partial t}+\left(\frac{\partial u}{\partial t} \right)^2-2u \frac{\partial^2 u}{\partial t^2},\\
   \frac{\partial^3 u}{\partial t^3}=&\frac{\partial u}{\partial S}+6u^2\frac{\partial u}{\partial t},
   \end{aligned}
  \right. 
\end{equation}
where the first equation in \eqref{eq:G(1,4)5} coincides with Chazy XI$(N=3)$ equation (see \cite{11}).

Adding each system in \eqref{eq:G(1,4)5}, we can obtain
\begin{equation}\label{eq:G(1,4)6}
   \frac{\partial^3 u}{\partial t^3} =\frac{3}{2}u^4+6u^2 \frac{\partial u}{\partial t}+\frac{1}{2}\left(\frac{\partial u}{\partial t} \right)^2-u \frac{\partial^2 u}{\partial t^2}+\frac{1}{2} \frac{\partial u}{\partial S}.
\end{equation}

\begin{question}
It is still an open question whether the equation \eqref{eq:G(1,4)6} coincides with which of the soliton equations.
\end{question}

\section{A generalization of the system \eqref{SdeG4} to three variables}
In this section, we present a generalization of the system \eqref{SdeG4} to three variables, which is equivalent to the polynomial Hamiltonian system
\begin{align}\label{SdeGaK}
\begin{split}
&dx=\frac{\partial H_1}{\partial y}dt+\frac{\partial H_2}{\partial y}ds+\frac{\partial H_3}{\partial y}du, \quad dy=-\frac{\partial H_1}{\partial x}dt-\frac{\partial H_2}{\partial x}ds-\frac{\partial H_3}{\partial x}du,\\
&dz=\frac{\partial H_1}{\partial w}dt+\frac{\partial H_2}{\partial w}ds+\frac{\partial H_3}{\partial w}du, \quad dw=-\frac{\partial H_1}{\partial z}dt-\frac{\partial H_2}{\partial z}ds-\frac{\partial H_3}{\partial z}du,\\
&dq=\frac{\partial H_1}{\partial p}dt+\frac{\partial H_2}{\partial p}ds+\frac{\partial H_3}{\partial p}du, \quad dp=-\frac{\partial H_1}{\partial q}dt-\frac{\partial H_2}{\partial q}ds-\frac{\partial H_3}{\partial q}du\\
\end{split}
\end{align}
with the symmetric Hamiltonians $H_i \in {\Bbb C}(t,s,u)[x,y,z,w,q,p] \ (i=1,2,3)$
\begin{align}\label{SdeGaKH}
\begin{split}
H_1 &=H_{II}(x,y,t;\alpha_3)+R(x,y,z,w,t,s;\alpha_3,\alpha_1)+R(x,y,q,p,t,u;\alpha_3,\alpha_4)\\
H_2 &=\pi(H_1), \quad H_3=(\pi \circ \pi)(H_1) \quad (\alpha_1+\alpha_2+\alpha_3+\alpha_4=1),
\end{split}
\end{align}
where the transformation $\pi$ is explicitly given by
\begin{equation}
\pi:(*) \rightarrow (z,w,q,p,x,y,s,u,t;\alpha_3,\alpha_2,\alpha_4,\alpha_1).
\end{equation}
Here the symbol $(*)$ denotes $(*):=(x,y,z,w,q,p,t,s,u;\alpha_1,\alpha_2,\alpha_3,\alpha_4)$, and the symbol\\
 $R(q_l,p_l,q_m,p_m,t_l,t_m;\alpha,\beta)$ is explicitly given by
\begin{align*}
&R(q_l,p_l,q_m,p_m,t_l,t_m;\alpha,\beta)\\
&=\frac{\beta}{t_l-t_m}q_lp_l-\frac{\beta}{t_l-t_m}p_lq_m-\frac{\alpha}{t_l-t_m}q_lp_m+\frac{\alpha}{t_l-t_m}q_mp_m-\frac{\{2(q_l-q_m)^2-t_l+t_m\}p_lp_m}{2(t_l-t_m)}.
\end{align*}

\begin{theorem}
Let us consider a polynomial Hamiltonian system with Hamiltonian\\
$H_i \in {\Bbb C}(t,s,u)[x,y,z,w,q,p] \ (i=1,2,3)$. We assume that

$(A1)$ $deg(H_i)=5$ with respect to $x,y,z,w,q,p$.

$(A2)$ This system becomes again a polynomial Hamiltonian system in each coordinate system $\{U_j,(x_j,y_j,z_j,w_j,q_j,p_j)\} \ (j=1,\dots,4)${\rm : \rm}
$$
{U_j}={{\Bbb C}^6} \ni {(x_j,y_j,z_j,w_j,q_j,p_j)} \ (j=1,\dots,4),
$$
via the following birational and symplectic transformations
\begin{align*}
&1)x_1=\frac{1}{x}, \quad y_1=-(yx+\alpha_3)x, \quad z_1=z, \quad w_1=w, \quad q_1=q, \quad p_1=p,\\
&2)x_2=\frac{1}{x}, \quad y_2=-\{(y+w+p-2x^2-t)x-2\left((z-x)x-\frac{t-s}{4}\right)\left(\frac{w}{x}\right)\\
&-2\left((q-x)x-\frac{t-u}{4}\right)\left(\frac{p}{x}\right)+\alpha_2\}x,\\
&z_2=\left((z-x)x-\frac{t-s}{2}\right)x, \quad w_2=\frac{w}{x^2}, \quad q_2=\left((q-x)x-\frac{t-u}{2}\right)x, \quad p_2=\frac{p}{x^2},\\
&3)x_3=x,\quad y_3=y,\quad z_3=\frac{1}{z},\quad w_3=-(zw+\alpha_1)z, \quad q_3=q, \quad p_3=p,\\
&4)x_4=x, \quad y_4=y, \quad z_4=z, \quad w_4=w, \quad q_4=\frac{1}{q}, \quad p_4=-(qp+\alpha_4)q.
\end{align*}
Then such a system coincides with the system \eqref{SdeGaK}.
\end{theorem}

\begin{theorem}
The system \eqref{SdeGaK} is invariant under the following transformations\rm{:\rm} with the notation $(*)=(x,y,z,w,q,p,t,s,u;\alpha_0,\alpha_1,\ldots,\alpha_4),$
\begin{align*}
        s_1: (*) &\rightarrow \left(x+\frac{\alpha_3}{y},y,z,w,q,p,t,s,u;\alpha_1,\alpha_2+\alpha_3,-\alpha_3,\alpha_4 \right),\\
        s_3: (*) &\rightarrow \left(x,y,z+\frac{\alpha_1}{w},w,q,p,t,s,u;-\alpha_1,\alpha_2+\alpha_1,\alpha_3,\alpha_4 \right),\\
        s_4: (*) &\rightarrow \left(x,y,z,w,q+\frac{\alpha_4}{p},p,t,s,u;\alpha_1,\alpha_2+\alpha_4,\alpha_3,-\alpha_4 \right),\\
        \pi_1: (*) &\rightarrow (z,w,q,p,x,y,s,u,t;\alpha_3,\alpha_2,\alpha_4,\alpha_1).
        \end{align*}
\end{theorem}

\section{Autonomous version of the degenerate Garnier system of type G(2,3) in two variables}

In this section, we find an autonomous version of the degenerate Garnier system of type G(2,3) in two variables given by
\begin{equation}\label{eq:AAA1}
  \left\{
  \begin{aligned}
   dq_1 =&\frac{\partial K_1}{\partial p_1}dt+\frac{\partial K_2}{\partial p_1}ds,\\
   dp_1 =&-\frac{\partial K_1}{\partial q_1}dt-\frac{\partial K_2}{\partial q_1}ds,\\
   dq_2 =&\frac{\partial K_1}{\partial p_2}dt+\frac{\partial K_2}{\partial p_2}ds,\\
   dp_2 =&-\frac{\partial K_1}{\partial q_2}dt-\frac{\partial K_2}{\partial q_2}ds
   \end{aligned}
  \right. 
\end{equation}
with the polynomial Hamiltonians
\begin{align}\label{eq:AAA2}
\begin{split}
K_1=&-2\eta p_1+2q_2^2 p_2^2+q_2^2 p_2+2(\alpha_0+\alpha_2)q_2 p_2+\alpha_0 q_2-2\eta q_1p_2+q_1p_1q_2-2p_1^2q_2,\\
K_2=&q_1^2 p_1-2q_1 p_1^2+\alpha_0 q_1-2(\alpha_0+\alpha_2)p_1-2\eta p_2-p_1q_2+q_1q_2p_2-4p_1q_2p_2.
\end{split}
\end{align}

\begin{proposition}
The system \eqref{eq:AAA1} satisfies the compatibility conditions$:$
\begin{equation}
\frac{\partial }{\partial s} \frac{\partial q_1}{\partial t}=\frac{\partial }{\partial t} \frac{\partial q_1}{\partial s}, \quad \frac{\partial }{\partial s} \frac{\partial p_1}{\partial t}=\frac{\partial }{\partial t} \frac{\partial p_1}{\partial s}, \quad \frac{\partial }{\partial s} \frac{\partial q_2}{\partial t}=\frac{\partial }{\partial t} \frac{\partial q_2}{\partial s}, \quad \frac{\partial }{\partial s} \frac{\partial p_2}{\partial t}=\frac{\partial }{\partial t} \frac{\partial p_2}{\partial s}.
\end{equation}
\end{proposition}

\begin{proposition}
The system \eqref{eq:AAA1} has $K_1$ and $K_2$ as its first integrals.
\end{proposition}

\begin{proposition}
Two Hamiltonians $K_1$ and $K_2$ satisfy
\begin{equation}
\{K_1,K_2\}=0,
\end{equation}
where
\begin{equation}
\{K_1,K_2\}=\frac{\partial K_1}{\partial p_1}\frac{\partial K_2}{\partial q_1}-\frac{\partial K_1}{\partial q_1}\frac{\partial K_2}{\partial p_1}+\frac{\partial K_1}{\partial p_2}\frac{\partial K_2}{\partial q_2}-\frac{\partial K_1}{\partial q_2}\frac{\partial K_2}{\partial p_2}.
\end{equation}
\end{proposition}
Here, $\{,\}$ denotes the poisson bracket such that $\{p_i,q_j\}={\delta}_{ij}$ (${\delta}_{ij}$:kronecker's delta).

\begin{theorem}\label{th:AAA1}
The system \eqref{eq:AAA1} is invariant under the birational transformations $s_0,s_1,s_2,s_3$ and $s_4$ with {\it the notation} $(*):=(q_1,p_1,q_2,p_2,\eta,t,s;\alpha_0,\alpha_1,\alpha_2)$\rm{; \rm}
\begin{align*}
s_0:(*) \rightarrow &\left(q_1,p_1+\frac{\eta}{q_2},q_2,p_2-\frac{\alpha_1}{q_2}-\frac{\eta q_1}{q_2^2},-\eta,t,s;\alpha_0+\alpha_1,-\alpha_1,\alpha_2+\alpha_1 \right),\\
s_1:(*) \rightarrow &\left(-q_1,p_1-\frac{1}{2}q_1,q_2,-p_2-\frac{1}{2},\eta,-t,s;-\alpha_2,-\alpha_1,-\alpha_0 \right),\\
s_2:(*) \rightarrow &\left(-\sqrt{-1}q_1,-\frac{1}{2}\sqrt{-1}(q_1-2p_1),-q_2,-p_2-\frac{1}{2},-\sqrt{-1}\eta,t,-\sqrt{-1}s;\alpha_2,\alpha_1,\alpha_0 \right),\\
s_3:(*) \rightarrow &\left(-q_1,-p_1,q_2,p_2,-\eta,t,-s;\alpha_0,\alpha_1,\alpha_2 \right),\\
s_4:(*) \rightarrow &\left(\sqrt{-1}q_1,\sqrt{-1}p_1,-q_2,p_2,-\sqrt{-1}\eta,-t,-\sqrt{-1}s;-\alpha_0,-\alpha_1,-\alpha_2 \right).
\end{align*}
\end{theorem}
Here, the parameters $\alpha_i$ satisfy the relation $\alpha_0+\alpha_1+\alpha_2=0$.

\begin{corollary}
The transformations described in Theorem \ref{th:AAA1} satisfy the following relations{\rm : \rm}
$$
s_0^2=s_1^2=s_3^2=1, \quad s_2^4=s_4^4=1.
$$
\end{corollary}

\begin{theorem}
Let us consider a polynomial Hamiltonian system with Hamiltonian $K \in {\Bbb C}[q_1,p_1,q_2,p_2]$. We assume that

$(C1)$ $deg(K)=6$ with respect to $q_1,p_1,q_2,p_2$.

$(C2)$ This system becomes again a polynomial Hamiltonian system in each coordinate $R_i \ (i=0,1,2)${\rm : \rm}
\begin{align*}
\begin{split}
&R_0:(x_0,y_0,z_0,w_0)=\left(\frac{1}{q_1},-(q_1p_1+q_2p_2+\alpha_0)q_1,\frac{q_2}{q_1},p_2q_1 \right),\\
&R_1:(x_1,y_1,z_1,w_1)=\left(q_1,p_1+\frac{\eta}{q_2},q_2,p_2-\frac{\alpha_1}{q_2}-\frac{\eta q_1}{q_2^2} \right),\\
&R_2:(x_2,y_2,z_2,w_2)=\left(-\left((q_1-2p_1)p_1-q_2 \left(p_2+\frac{1}{2} \right)-\alpha_2 \right)p_1,\frac{1}{p_1},\frac{q_2}{p_1},\left(p_2+\frac{1}{2} \right)p_1 \right),
\end{split}
\end{align*}
where the parameters $\alpha_i$ satisfy the relation $\alpha_0+\alpha_1+\alpha_2=0$. Then such a system coincides with the Hamiltonian system \eqref{eq:AAA1} with the polynomial Hamiltonians $K_1,K_2$.
\end{theorem}
We remark that under the above assumptions $(C1)$ and $(C2)$ we can obtain three polynomial Hamiltonians $K_1,K_2$ and $K_3$. We easily see that
$$
K_3=K_2^2-2(\alpha_1+\alpha_2) K_1.
$$
We note that the conditions $(C2)$ should be read that
\begin{align*}
\begin{split}
&R_i(K) \quad (i=0,1,2)
\end{split}
\end{align*}
are polynomials with respect to $x_i,y_i,z_i,w_i$.

These holomorphy conditions $R_0,R_1$ and $R_2$ are new.

\section{Autonomous version of the degenerate Garnier system of type G(1,4) in two variables}

In this section, we find an autonomous version of the degenerate Garnier system of type G(1,4) in two variables given by
\begin{equation}\label{eq:A1}
  \left\{
  \begin{aligned}
   dq_1 =&\frac{\partial K_1}{\partial p_1}dt+\frac{\partial K_2}{\partial p_1}ds,\\
   dp_1 =&-\frac{\partial K_1}{\partial q_1}dt-\frac{\partial K_2}{\partial q_1}ds,\\
   dq_2 =&\frac{\partial K_1}{\partial p_2}dt+\frac{\partial K_2}{\partial p_2}ds,\\
   dp_2 =&-\frac{\partial K_1}{\partial q_2}dt-\frac{\partial K_2}{\partial q_2}ds
   \end{aligned}
  \right. 
\end{equation}
with the polynomial Hamiltonians
\begin{align}\label{eq:A2}
\begin{split}
K_1=&-q_1^2p_1+p_1^2-\alpha_0 q_1-g p_1-q_2 p_2^2+\alpha_1 p_2+p_1 q_2-q_1q_2p_2,\\
K_2=&-\alpha_1 p_1+q_2^2 p_2-g q_2 p_2+\alpha_0 q_2-\alpha_1 q_1 p_2+q_1 p_1q_2+2p_1q_2p_2+q_1q_2p_2^2 \quad (g \in {\Bbb C}).
\end{split}
\end{align}

\begin{proposition}
The system \eqref{eq:A1} satisfies the compatibility conditions$:$
\begin{equation}
\frac{\partial }{\partial s} \frac{\partial q_1}{\partial t}=\frac{\partial }{\partial t} \frac{\partial q_1}{\partial s}, \quad \frac{\partial }{\partial s} \frac{\partial p_1}{\partial t}=\frac{\partial }{\partial t} \frac{\partial p_1}{\partial s}, \quad \frac{\partial }{\partial s} \frac{\partial q_2}{\partial t}=\frac{\partial }{\partial t} \frac{\partial q_2}{\partial s}, \quad \frac{\partial }{\partial s} \frac{\partial p_2}{\partial t}=\frac{\partial }{\partial t} \frac{\partial p_2}{\partial s}.
\end{equation}
\end{proposition}

\begin{proposition}
The system \eqref{eq:A1} has $K_1$ and $K_2$ as its first integrals.
\end{proposition}

\begin{proposition}
Two Hamiltonians $K_1$ and $K_2$ satisfy
\begin{equation}
\{K_1,K_2\}=0,
\end{equation}
where
\begin{equation}
\{K_1,K_2\}=\frac{\partial K_1}{\partial p_1}\frac{\partial K_2}{\partial q_1}-\frac{\partial K_1}{\partial q_1}\frac{\partial K_2}{\partial p_1}+\frac{\partial K_1}{\partial p_2}\frac{\partial K_2}{\partial q_2}-\frac{\partial K_1}{\partial q_2}\frac{\partial K_2}{\partial p_2}.
\end{equation}
\end{proposition}
Here, $\{,\}$ denotes the poisson bracket such that $\{p_i,q_j\}={\delta}_{ij}$ (${\delta}_{ij}$:kronecker's delta).

\begin{theorem}\label{th:A1}
The system \eqref{eq:A1} is invariant under the birational and symplectic transformations $s_0,s_1,s_2,\pi:$ with {\it the notation} $(*):=(q_1,p_1,q_2,p_2,t,s;\alpha_0,\alpha_1,\alpha_2)$\rm{; \rm}
\begin{align*}
s_0:(*) \rightarrow &\left(q_1+\frac{q_2p_2+\alpha_0}{p_1},p_1,-\frac{p_2(q_2p_2-\alpha_1)}{p_1},\frac{p_1}{p_2},t,-s;-\alpha_0-\alpha_1,\alpha_1,\alpha_0 \right),\\
s_1:(*) \rightarrow &\left(q_1,p_1,q_2,p_2-\frac{\alpha_1}{q_2},t,s;\alpha_0+\alpha_1,-\alpha_1,\alpha_2 \right),\\
s_2:(*) \rightarrow &(q_1+\frac{(q_1+p_2)q_2+\alpha_2}{p_1+q_2-q_1^2-g},\\
&p_1+q_2-q_1^2-\frac{(q_1+p_2)\{(q_1+p_2)q_2-\alpha_1\}}{p_1+q_2-q_1^2-g}+\left(q_1+\frac{(q_1+p_2)q_2+\alpha_2}{p_1+q_2-q_1^2-g} \right)^2,\\
&\frac{(q_1+p_2)\{(q_1+p_2)q_2-\alpha_1\}}{p_1+q_2-q_1^2-g},-q_1-\frac{p_1+q_2-q_1^2-g}{q_1+p_2}-\frac{(q_1+p_2)q_2+\alpha_2}{p_1+q_2-q_1^2-g},t,-s;\\
&\alpha_2,\alpha_1,-\alpha_2-\alpha_1),\\
\pi:(*) \rightarrow &(-q_1,-(p_1+q_2-q_1^2-g),q_2,p_2+q_1,t,-s;\alpha_2,\alpha_1,\alpha_0).
\end{align*}
\end{theorem}
Here, the parameters $\alpha_i$ satisfy the relation $\alpha_0+\alpha_1+\alpha_2=0$.

\begin{corollary}
The transformations described in Theorem \ref{th:A1} satisfy the following relations{\rm : \rm}
$$
s_0^2=s_1^2=s_2^2=\pi^2=(s_1s_0)^2=1, \quad s_0 \pi=\pi s_2.
$$
\end{corollary}

\begin{proposition}
The transformation $s_0$ can be obtained by composing the following transformations$:$

{\bf Step 1:} We transform the system \eqref{th:A1} by the birational and symplectic transformation$:$
$$
g_1:(x_1,y_1,z_1,w_1)=\left(q_1,p_1,-(q_2p_2-\alpha_1)p_2,\frac{1}{p_2} \right).
$$

{\bf Step 2:} We then transform the system obtained by Step 1 by the birational and symplectic transformation$:$
$$
g_2:(x_2,y_2,z_2,w_2)=\left(x_1-\frac{z_1w_1-\alpha_0-\alpha_1}{y_1},y_1,\frac{z_1}{y_1},w_1y_1 \right).
$$
Thus, we can obtain the B\"acklund transformation $s_0$ of the system \eqref{th:A1}.
\end{proposition}

\begin{theorem}
Let us consider a polynomial Hamiltonian system with Hamiltonian $K \in {\Bbb C}[q_1,p_1,q_2,p_2]$. We assume that

$(C1)$ $deg(K)=5$ with respect to $q_1,p_1,q_2,p_2$.

$(C2)$ This system becomes again a polynomial Hamiltonian system in each coordinate $R_i \ (i=0,1,2)${\rm : \rm}
\begin{align*}
\begin{split}
R_0:(x_0,y_0,z_0,w_0)=&\left(\frac{1}{q_1},-(q_1p_1+q_2p_2+\alpha_0)q_1,\frac{q_2}{q_1},p_2q_1 \right),\\
R_1:(x_1,y_1,z_1,w_1)=&\left(q_1,p_1,-(q_2p_2-\alpha_1)p_2,\frac{1}{p_2} \right),\\
R_2:(x_2,y_2,z_2,w_2)=&\left(\frac{1}{q_1},-(q_1(p_1+q_2-q_1^2-g)+q_2(q_1+p_2)+\alpha_2)q_1,\frac{q_2}{q_1},(p_2+q_1)q_1 \right),
\end{split}
\end{align*}
where the parameters $\alpha_i$ satisfy the relation $\alpha_0+\alpha_1+\alpha_2=0$. Then such a system coincides with the Hamiltonian system \eqref{eq:A1} with the polynomial Hamiltonians $K_1,K_2$.
\end{theorem}
We note that the conditions $(C2)$ should be read that
\begin{align*}
\begin{split}
&R_0(K), \quad R_1(K), \quad R_2(K)
\end{split}
\end{align*}
are polynomials with respect to $x_i,y_i,z_i,w_i$.

\begin{theorem}
The birational transformations
\begin{equation}\label{eq:A3}
  \left\{
  \begin{aligned}
   x =&q_1,\\
   y =&q_2+2p_1-q_1^2-g,\\
   z =&2g q_1+2q_1^3-3q_1q_2+2\alpha_0+\alpha_1,\\
   w =&-2g^2-8gq_1^2-6q_1^4+4gp_1+12q_1^2 p_1+5gq_2+6q_1q_2p_2+12q_1^2q_2-6p_1q_2-3q_2^2-3\alpha_1 q_1
   \end{aligned}
  \right. 
\end{equation}
take the Hamiltonian system \eqref{eq:A1} to the system
\begin{equation}\label{eq:A4}
  \left\{
  \begin{aligned}
   dx =&y dt+\left(-\frac{4}{3}xy+\frac{w}{3x}-\frac{y(z-2\alpha_0-\alpha_1)}{3x^2} \right)ds,\\
   dy =&z dt+f_1(x,y,z,w)ds,\\
   dz =&w dt+f_2(x,y,z,w)ds,\\
   dw =&\frac{1}{3}g(2\alpha_0+\alpha_1)-\frac{2}{3}g^2x-\frac{1}{3}gz-\frac{28}{3}\alpha_0 x^2+\frac{13}{3}\alpha_1 x^2-\frac{34}{3}gx^3-\frac{32}{3}x^5+9xy^2+\frac{26}{3}x^2z\\
&+\frac{(2\alpha_0+\alpha_1)^2}{3x}-\frac{5(2\alpha_0+\alpha_1)z}{3x}+\frac{yw}{x}+\frac{4z^2}{3x}+\frac{gy^2}{x}+\frac{3(2\alpha_0+\alpha_1)y^2}{2x^2}-\frac{3y^2 z}{2x^2}\\
&-\frac{(w-6gx^2-6x^4-2gy-6x^2 y+3xz-6\alpha_0 x)}{2(2gx+2x^3-z+2\alpha_0+\alpha_1)} \times\\
&(w+6gx^2+6x^4-2gy-6x^2 y-3xz+6\alpha_0 x) dt+f_3(x,y,z,w)ds,
   \end{aligned}
  \right. 
\end{equation}
where $f_i(x,y,z,w) \in {\Bbb C}[x,y,z,w] \ (i=1,2,3)$. 
\end{theorem}
Setting $u:=x$, we see that
\begin{equation}
\frac{\partial u}{\partial t}=y, \quad  \frac{\partial^2 u}{\partial t^2}=z, \quad \frac{\partial^3 u}{\partial t^3}=w,
\end{equation}
and
\begin{equation}\label{eq:A5}
  \left\{
  \begin{aligned}
   \frac{\partial^4 u}{\partial t^4} =&\frac{1}{3}g(2\alpha_0+\alpha_1)-\frac{2}{3}g^2u-\frac{1}{3}g\frac{\partial^2 u}{\partial t^2}-\frac{28}{3}\alpha_0 u^2+\frac{13}{3}\alpha_1 u^2-\frac{34}{3}gu^3-\frac{32}{3}u^5+9u \left(\frac{\partial u}{\partial t} \right)^2\\
&+\frac{26}{3}u^2\frac{\partial^2 u}{\partial t^2}+\frac{(2\alpha_0+\alpha_1)^2}{3u}-\frac{5(2\alpha_0+\alpha_1)}{3u}\frac{\partial^2 u}{\partial t^2}+\frac{1}{u}\frac{\partial u}{\partial t} \frac{\partial^3 u}{\partial t^3}+\frac{4}{3u} \left(\frac{\partial^2 u}{\partial t^2} \right)^2\\
&+\frac{g}{u} \left(\frac{\partial u}{\partial t} \right)^2+\frac{3(2\alpha_0+\alpha_1)}{2u^2} \left(\frac{\partial u}{\partial t} \right)^2-\frac{3}{2u^2} \left(\frac{\partial u}{\partial t} \right)^2 \frac{\partial^2 u}{\partial t^2}\\
&-\frac{(\frac{\partial^3 u}{\partial t^3}-6gu^2-6u^4-2g\frac{\partial u}{\partial t}-6u^2 \frac{\partial u}{\partial t}+3u\frac{\partial^2 u}{\partial t^2}-6\alpha_0 u)}{2(2gu+2u^3-\frac{\partial^2 u}{\partial t^2}+2\alpha_0+\alpha_1)} \times\\
&\left(\frac{\partial^3 u}{\partial t^3}+6gu^2+6u^4-2g\frac{\partial u}{\partial t}-6u^2 \frac{\partial u}{\partial t}-3u\frac{\partial^2 u}{\partial t^2}+6\alpha_0 u \right),\\
   \frac{\partial u}{\partial s} =&-\frac{4}{3}u \frac{\partial u}{\partial t}+\frac{1}{3u}\frac{\partial^3 u}{\partial t^3}-\frac{1}{3u^2}\frac{\partial u}{\partial t} \left(\frac{\partial^2 u}{\partial t^2}-2\alpha_0-\alpha_1 \right).
   \end{aligned}
  \right. 
\end{equation}
The first equation in \eqref{eq:A5} is the fourth-order ordinary differential equation of rational type. This equation is new. We see that the second equation in \eqref{eq:A5} can be considered as homogeneous polynomial of degree 4 when we set  $[\alpha_0]=[\alpha_1]=3$, $[u]=1,[\frac{\partial u}{\partial t}]=2,[\frac{\partial^2 u}{\partial t^2}]=3,[\frac{\partial^3 u}{\partial t^3}]=4$ and $[\frac{\partial u}{\partial s}]=3$.

For the second equation in \eqref{eq:A5}, we see that this system admits travelling wave solutions $u(t,s)=U(t-cs)$, where $U(T)$ satisfies the equation
\begin{equation}\label{eq:A6}
\frac{d^2 U}{d T^2}=2U^3-3cU^2+2\alpha_0+\alpha_1.
\end{equation}
This equation is an autonomous version of the second Painlev\'e equation (see \cite{Ince}).

\begin{question}
It is still an open question whether the second equation in \eqref{eq:A5} coincides with which of the soliton equations.
\end{question}

\section{Autonomous version of the degenerate Garnier system of type G(5) in two variables}

In this section, we find an autonomous version of the degenerate Garnier system of type G(5) in two variables given by
\begin{equation}\label{eq:AA1}
  \left\{
  \begin{aligned}
   dq_1 =&\frac{\partial K_1}{\partial p_1}dt+\frac{\partial K_2}{\partial p_1}ds,\\
   dp_1 =&-\frac{\partial K_1}{\partial q_1}dt-\frac{\partial K_2}{\partial q_1}ds,\\
   dq_2 =&\frac{\partial K_1}{\partial p_2}dt+\frac{\partial K_2}{\partial p_2}ds,\\
   dp_2 =&-\frac{\partial K_1}{\partial q_2}dt-\frac{\partial K_2}{\partial q_2}ds
   \end{aligned}
  \right. 
\end{equation}
with the polynomial Hamiltonians
\begin{align}\label{eq:AA2}
\begin{split}
K_1=&q_2^2 p_2+\alpha_0 q_2-q_1 p_2+p_1p_2+q_1 p_1 q_2+\frac{1}{2}p_1^2 q_2,\\
K_2=&q_1^2 p_1-\frac{1}{2}q_1 p_1^2+\alpha_0 q_1+\frac{1}{2}p_2^2+q_1q_2p_2+p_1q_2p_2+\frac{1}{2}p_1^2 q_2^2.
\end{split}
\end{align}

\begin{proposition}
The system \eqref{eq:AA1} satisfies the compatibility conditions$:$
\begin{equation}
\frac{\partial }{\partial s} \frac{\partial q_1}{\partial t}=\frac{\partial }{\partial t} \frac{\partial q_1}{\partial s}, \quad \frac{\partial }{\partial s} \frac{\partial p_1}{\partial t}=\frac{\partial }{\partial t} \frac{\partial p_1}{\partial s}, \quad \frac{\partial }{\partial s} \frac{\partial q_2}{\partial t}=\frac{\partial }{\partial t} \frac{\partial q_2}{\partial s}, \quad \frac{\partial }{\partial s} \frac{\partial p_2}{\partial t}=\frac{\partial }{\partial t} \frac{\partial p_2}{\partial s}.
\end{equation}
\end{proposition}

\begin{proposition}
The system \eqref{eq:AA1} has $K_1$ and $K_2$ as its first integrals.
\end{proposition}

\begin{proposition}
Two Hamiltonians $K_1$ and $K_2$ satisfy
\begin{equation}
\{K_1,K_2\}=0,
\end{equation}
where
\begin{equation}
\{K_1,K_2\}=\frac{\partial K_1}{\partial p_1}\frac{\partial K_2}{\partial q_1}-\frac{\partial K_1}{\partial q_1}\frac{\partial K_2}{\partial p_1}+\frac{\partial K_1}{\partial p_2}\frac{\partial K_2}{\partial q_2}-\frac{\partial K_1}{\partial q_2}\frac{\partial K_2}{\partial p_2}.
\end{equation}
\end{proposition}
Here, $\{,\}$ denotes the poisson bracket such that $\{p_i,q_j\}={\delta}_{ij}$ (${\delta}_{ij}$:kronecker's delta).

\begin{theorem}\label{th:AA1}
The system \eqref{eq:AA1} is invariant under the birational transformations $s_0,s_1$ with {\it the notation} $(*):=(q_1,p_1,q_2,p_2,t,s;\alpha_0,\alpha_1)$\rm{; \rm}
\begin{align*}
s_0:(*) \rightarrow &(q_1,-(p_1-2q_1+2q_2^2),-q_2,p_2+4q_1q_2-2q_2^3,t,-s;\alpha_0,\alpha_1 ),\\
s_1:(*) \rightarrow &\left(\sqrt{-1}q_1,\frac{p_1-2q_1+2q_2^2}{\sqrt{-1}},(-1)^{\frac{1}{4}}q_2,\frac{p_2+4q_1q_2-2q_2^3}{(-1)^{\frac{1}{4}}},(-1)^{\frac{3}{4}}t,\sqrt{-1}s;-\alpha_0,-\alpha_1 \right).
\end{align*}
\end{theorem}
Here, the parameters $\alpha_i$ satisfy the relation $\alpha_0+\alpha_1=0$.

\begin{corollary}
The transformations described in Theorem \ref{th:AA1} satisfy the following relations{\rm : \rm}
$$
s_0^2=1, \quad s_1^8=1.
$$
\end{corollary}

\begin{theorem}
Let us consider a polynomial Hamiltonian system with Hamiltonian $K \in {\Bbb C}[q_1,p_1,q_2,p_2]$. We assume that

$(C1)$ $deg(K)=6$ with respect to $q_1,p_1,q_2,p_2$.

$(C2)$ This system becomes again a polynomial Hamiltonian system in each coordinate $R_i \ (i=0,1)${\rm : \rm}
\begin{align*}
\begin{split}
&R_0:(x_0,y_0,z_0,w_0)=\left(\frac{q_1}{q_2},p_1q_2,\frac{1}{q_2},-(q_2p_2+q_1p_1+\alpha_0)q_2 \right),\\
&R_1:(x_1,y_1,z_1,w_1)\\
&=\left(\frac{q_1}{q_2},(p_1-2q_1+2q_2^2)q_2,\frac{1}{q_2},-((p_2+4q_1q_2-2q_2^3)q_2+q_1(p_1-2q_1+2q_2^2)+\alpha_1)q_2 \right),
\end{split}
\end{align*}
where the parameters $\alpha_i$ satisfy the relation $\alpha_0+\alpha_1=0$. Then such a system coincides with the Hamiltonian system \eqref{eq:AA1} with the polynomial Hamiltonians $K_1,K_2$.
\end{theorem}
We remark that under the above assumptions $(C1)$ and $(C2)$ we can obtain three polynomial Hamiltonians $K_1,K_2$ and $K_3$. We easily see that
$$
K_3=K_1^2-2\alpha_0 K_2.
$$
We note that the conditions $(C2)$ should be read that
\begin{align*}
\begin{split}
&R_0(K), \quad R_1(K)
\end{split}
\end{align*}
are polynomials with respect to $x_i,y_i,z_i,w_i$.

These holomorphy conditions $R_0$ and $R_1$ are new.

\begin{theorem}
The birational transformations
\begin{equation}\label{eq:AA3}
  \left\{
  \begin{aligned}
   x =&-(3q_1-2q_2^2)q_2,\\
   y =&-q_1+p_1+q_2^2,\\
   z =&q_2,\\
   w =&-3(-q_1^2+q_1p_1+q_2p_2+4q_1q_2^2-p_1q_2^2-2q_2^4)
   \end{aligned}
  \right. 
\end{equation}
take the Hamiltonian system \eqref{eq:AA1} to the system
\begin{equation}\label{eq:AA4}
  \left\{
  \begin{aligned}
   dz =&y dt+\left(\frac{4}{3}yz-\frac{w}{3z}+\frac{1}{3z^2}xy \right)ds,\\
   dy =&x dt+f_1(x,y,z,w)ds,\\
   dx =&w dt+f_2(x,y,z,w)ds,\\
   dw =&\left(-\frac{2xy^2}{z^2}+\frac{3x^2}{2z}+\frac{2yw}{z}+\frac{5}{2}y^2z+5xz^2-\frac{5}{2}z^5+3\alpha_0 z \right)dt+f_3(x,y,z,w)ds,
   \end{aligned}
  \right. 
\end{equation}
where $f_i(x,y,z,w) \in {\Bbb C}[x,y,z,w] \ (i=1,2,3)$. 
\end{theorem}
Setting $u:=z$, we see that
\begin{equation}
\frac{\partial u}{\partial t}=y, \quad  \frac{\partial^2 u}{\partial t^2}=x, \quad \frac{\partial^3 u}{\partial t^3}=w,
\end{equation}
and
\begin{equation}\label{eq:AA5}
  \left\{
  \begin{aligned}
   \frac{\partial^4 u}{\partial t^4} =&-\frac{2}{u^2}\left(\frac{\partial u}{\partial t} \right)^2 \frac{\partial^2 u}{\partial t^2}+\frac{3}{2u} \left(\frac{\partial^2 u}{\partial t^2} \right)^2+\frac{2}{u}\frac{\partial u}{\partial t} \frac{\partial^3 u}{\partial t^3} +\frac{5}{2} \left(\frac{\partial u}{\partial t} \right)^2u+5\frac{\partial^2 u}{\partial t^2} u^2\\
&-\frac{5}{2}u^5+3\alpha_0 u,\\
   \frac{\partial u}{\partial s} =&\frac{4}{3}u\frac{\partial u}{\partial t}-\frac{1}{3u}\frac{\partial^3 u}{\partial t^3}+\frac{1}{3u^2}\frac{\partial u}{\partial t}\frac{\partial^2 u}{\partial t^2}.
   \end{aligned}
  \right. 
\end{equation}
The first equation in \eqref{eq:AA5} is the fourth-order ordinary differential equation of rational type. This equation is new (see \cite{11}).

\end{document}